\documentclass[reqno]{amsart}
\usepackage{graphicx,xcolor,cite,mathtools,bbold,bbding}
\usepackage{amssymb,amsmath,amsthm,mathrsfs,paralist,bm,esint,setspace}
\usepackage[ngerman,french,english]{babel}
\RequirePackage[colorlinks,citecolor=blue,urlcolor=blue]{hyperref}
\usepackage{cite}
\usepackage{cases}
\usepackage{tcolorbox}
\usepackage{tikz}
\usepackage{pgfplots}
\usepackage{subcaption}
\usetikzlibrary{arrows.meta}
\pgfplotsset{compat=1.14}
\pgfplotsset{every x tick label/.append style={font=\footnotesize, yshift=0.6ex}}
\pgfplotsset{every y tick label/.append style={font=\footnotesize, xshift=0.5ex}}
\DeclareMathOperator{\Var}{Var}

\DeclareMathOperator{\Cov}{Cov}

\DeclareMathOperator{\dimh}{dim_{_{\rm H}}}
\DeclareMathOperator{\dimm}{\overline{dim}_{_{\rm M}}}

\DeclareMathOperator{\fS}{\mathfrak S}

\newcommand{\N}{\mathbb{N}}
\newcommand{\Z}{\mathbb{Z}}
\newcommand{\Q}{\mathbb{Q}}
\newcommand{\R}{\mathbb{R}}

\newcommand{\lip}{\text{Lip}}
\renewcommand{\P}{\mathrm{P}}

\newcommand{\E}{\mathrm{E}}

\newcommand{\1}{\mathbb{1}}
\renewcommand{\d}{{\rm d}}

\newcommand{\e}{{\rm e}}

\renewcommand{\ge}{\geqslant}
\renewcommand{\le}{\leqslant}
\definecolor{CYL}{rgb}{0.3,0.1,0.1}
\definecolor{DK}{rgb}{0.5,0.3,0.5}
\author[D. Khoshnevisan and C.Y. Lee]{Davar Khoshnevisan \and Cheuk Yin Lee}
\address{Department of Mathematics, The University of Utah, Salt Lake City, Utah 84112-0090,
	USA}
\email{davar@math.utah.edu}
\address{School of Science and Engineering, The Chinese University of Hong Kong (Shenzhen), Longgang,
	Shenzhen, Guangdong, 518172, China}
\email{leecheukyin@cuhk.edu.cn}
\title[Slow growth for SPDEs]{Points of slow growth for parabolic SPDEs}\thanks{
	Research supported in part by the National Science Foundation grant
        DMS-2245242 and the Shenzhen Peacock grant 2025TC0013.}
\newtheorem{stat}{Statement}[section]
\newtheorem{proposition}[stat]{Proposition}

\newtheorem{corollary}[stat]{Corollary}
\newtheorem{theorem}[stat]{Theorem}
\newtheorem{lemma}[stat]{Lemma}
\theoremstyle{definition}

\newtheorem{conjecture}[stat]{Conjecture}
\newtheorem{remark}[stat]{Remark}

\numberwithin{equation}{section}
\usepackage[cal=dutchcal,calscaled=1.05,
	scr=boondoxo,scrscaled=1.05,
	]{%
mathalfa}
\date{August 10, 2026}
\begin{document}
\maketitle
%\tableofcontents
%----------------------------------------------------------------------------------------------------
\setcounter{tocdepth}{3}% to get subsubsections in toc
\let\oldtocsection=\tocsection
\let\oldtocsubsection=\tocsubsection
\let\oldtocsubsubsection=\tocsubsubsection
%----------------------------------------------------------------------------------------------------

\begin{abstract}
	Consider the stochastic PDE, $\partial_tu = \partial^2_x u + \sigma(u) \dot{W}$
	on $\R_+\times\R$,
	subject to $u(0)\equiv1$, where $\dot{W}$ denotes space-time white noise
	on $\R_+\times\R$ and $\sigma:\R\to\R$ is Lipschitz continuous. It is known
	that $u(t\,,x)-1$ has approximately a Gaussian distribution for every $x$ 
	when $t\approx0$. Here we prove that there exist random points $x\in\R$
	where the fluctuations of the solution near times zero are almost surely
	of sharp order $t^{1/4}$. Our work 
	bears some loose resemblance to the study of the slow points
	of Brownian motion increments \cite{Dvoretzky,
	BassBurdzy,Davis1983,DavisPerkins,
	GreenwoodPerkins,Kahane1974,
	Kahane1976,Kahane1983,Perkins1983}, though significant challenges
	arise due to the infinite-dimensional nature of the present problem.
\end{abstract}

\section{Introduction}

Let us consider the stochastic partial differential equation (SPDE, for short),
\begin{equation}\label{u}\left[\begin{aligned}
	&\partial_t u(t\,,x) = \partial^2_x u(t\,,x) + \sigma(u(t\,,x)) \dot{W}(t\,,x) 
		&\qquad\text{for $(t\,,x)\in(0\,,\infty)\times\R$},\\
	&u(0) = 1&\text{on $\R$},
\end{aligned}\right.\end{equation}
where $\sigma:\R\to\R$ is a nonrandom Lipschitz continuous function,
and $\dot{W}$ denotes a space-time white noise. 
In order to avoid degeneracies, we will also assume that
$\sigma(1)\neq0$,
for  $u\equiv 1$ otherwise.

It is well known that  \eqref{u} is well posed
and has a random-field solution $u=\{u(t\,,x)\}_{t\ge0,x\in\R}$; 
see Dalang \cite{Dalang1999}, and that the solution $u$ is locally 
H\"older continuous in its variables, as well. In fact,
$u\in C^{a/2,a}_{\textit{loc}}(\R_+\times\R)$ a.s.~for every $a\in(0\,,\frac12)$; 
see the proof of Theorem 3.8 of Walsh \cite{Walsh} or 
Dalang and Sanz-Sol\'e \cite[Proposition 3.3.8]{DS}.
The continuity of $u$ implies that $u(t\,,x)\approx1$ to leading term,
when $t\approx0$, valid for every $x\in\R$. The error in this approximation
describes the rate of growth of the solution, at the spatial point $x$,
away from its initial profile. Also, that error is known to be of sharp order $t^{1/4}$.
For instance, it is known that
$[u(t\,,x)-1]/t^{1/4}$ converges in distribution to a non-degenerate normal
law for every $x\in\R$; see  Amir, Corwin, and Quastel \cite{ACQ},
Hairer and Pardoux \cite{HairerPardoux}, and Khoshnevisan, Swanson, 
Xiao, and Zhang \cite{KSXZ}.
In particular, it follows from this that
\begin{equation}\label{typical}
	\P\left\{\limsup_{t\to0^+}\frac{u(t\,,x)-1}{t^{1/4}} =\infty\right\}=1
	\quad\forall x\in\R.
\end{equation}
The above describes the growth of $u(t\,,x)$ away from its 
starting point $1$ at typical points $x\in\R$.
We say that a point $x\in\R$ is a \emph{point of slow growth}
for \eqref{u} if $|u(t\,,x)-1|=\mathcal{O}(t^{1/4})$ as $t\to0^+$
a.s. 
With these comments in mind, let us consider the random set
\begin{align}\label{S_u(theta)}
	\fS(\theta) = \left\{ x \in\R :\ \limsup_{t \to 0^+} 
	\frac{|u(t\,,x)-1|}{t^{1/4}} \le |\sigma(1)| \theta \right\},
\end{align}
where $\theta>0$ is nonrandom but otherwise arbitrary. 
We think of the elements of $\fS(\theta)$ -- if there are any --
as the \emph{$\theta$-slow points for \eqref{u}.}
Thanks to \eqref{typical} and Fubini's theorem, points of slow growth
are a.s.~Lebesgue-null. As we see next, they exist nevertheless,
and there are no points that grow more slowly than those in $\cup_{\theta\ge\theta_c}
\fS(\theta)$ for a certain number $\theta_c \in (0\,,\infty)$ that
does not depend on the nonlinearity $\sigma$;
see \eqref{theta_c} below. 
Let us first recall that ``$f(x) = \mathscr{o}(g(x))$ as $x\to a$'' 
means that $\lim_{x\to a}|f(x)/g(x)|=0$.
The following is the main
contribution of this paper.

\begin{theorem}\label{th:BC:u}
	If the underlying probability space $(\Omega\,,\mathscr{F},\P)$ is 
	complete, then:
	\begin{enumerate}[\noindent (1)]
	\item There exists a strictly decreasing, convex, 
		and continuous function $\lambda: (0\,,\infty) \to (0\,,\infty)$, 
		independent of $\sigma$, such that $\lim_{r\to0^+} \lambda(r) = \infty$, 
		$\lim_{r \to \infty} \lambda(r) = 0$, and 
		\[
			\P\left\{ |u(t\,,0)-1| \le |\sigma(1)| \theta t^{1/4}\ \ 
			\forall t \in [\varepsilon f(\varepsilon),\varepsilon]\right\} 
			= f(\varepsilon)^{\lambda(\theta)+\mathscr{o}(1)}\quad
			\text{as $\varepsilon\to0$},
		\]
		for every function $f:(0\,,1/\e)\to(0\,,1/\e)$ that vanishes at zero and satisfies
		$f(\varepsilon) \ge \varepsilon^{\ell+\mathscr{o}(1)}$, as $\varepsilon\to0^+$,
		for a fixed number $\ell=\ell(f)>0$.
	\item Let $K$ denote a nonrandom compact set in $\R$. Then,
		\[
			\P\{\fS(\theta)\cap K\neq\varnothing\}=
			\begin{cases}
				1 & \text{if $\lambda(\theta) < \frac12\dimh K$},\\
				0 & \text{if $\lambda(\theta) > \frac12 \dimm K$},
			\end{cases}
		\]
		where $\dimm$ and $\dimh$ denote respectively the upper Minkowski, or box, dimension
		and the Hausdorff dimension (see Falconer \textnormal{\cite{Falconer}}); and
	\item  There exists a single $\P$-null set off which $\dimh \fS(\theta) = 1 - 2\lambda(\theta)$
		for every 
		\begin{equation}\label{theta_c}
			\theta\ge \theta_c := \lambda^{-1}(\tfrac12).
		\end{equation}
	\end{enumerate}
\end{theorem}

Corresponding results for ``points of fast growth'' can be found in Huang and Khoshnevisan
\cite{HuangKh}.

\begin{remark}\label{rem:meas}
	As usual, the completeness condition for the underlying probability space
	can be made without loss in generality. It is used here in order to overcome 
	potential measurability issues. For instance, we establish Part (2) by proving 
	that $\{\omega\in\Omega:\ \fS(\theta)\cap K\neq\varnothing\}$
	is a subset of a $\P$-null set when $\lambda(\theta) > \frac12\dimm K$,
	and contains a full-probability set when $\lambda(\theta) < \frac12\dimh K$.
	Likewise, Part (3) is proved by showing
	that the set $\{\omega\in\Omega:\ \dimh \fS(\theta)(\omega) \neq 1 - 
	2\lambda(\theta)\}$ lies in a $\P$-null set. In the sequel, we 
	might suppress discussions
	on measurability. In all such cases, measurability issues can be verified individually using this type of
	inclusion/exclusion idea.
\end{remark}

Part (1) of this theorem is basically an analytic fact, and builds on our
recent work on Gaussian processes \cite{KL}. It also can be viewed as a contribution
to weighted small-ball problems for nonlinear SPDEs; compare with the recent works of Athreya, Joseph, and Mueller \cite{AthreyaJosephMueller},
Chen \cite{Chen},
Foondun, Joseph, and Kim \cite{FoondunJosephKim},
Guo, Song, Wang, and Xiao \cite{GuoSongWangXiao},
Hu and Lee \cite{HuLee},
Khoshnevisan, Kim, and Mueller \cite{KKimMueller},
and Martin \cite{Martin}.
In light of Part (1), we can consider $K$
to be an arbitrary compact interval in $\R$ in order to immediately deduce from
Part (2) of Theorem \ref{th:BC:u} two assertions about the sample
functions of the solution to \eqref{u}: 
\begin{enumerate}
\item[$(\mathscr{a})$]  Almost surely, $\theta$-slow points
	exist when $\theta>\theta_c$ [see \eqref{theta_c}]; and 
\item[$(\mathscr{b})$] Almost surely, $\theta$-slow points
	do not exist when $\theta<\theta_c$. 
\end{enumerate}

There are many analogies between Theorem \ref{th:BC:u}
and the findings of the literature on the slow points of
Brownian motion.
Item $(\mathscr{a})$ bears resemblance
to the well-known fact that Brownian increments can be as slow as
roughly $t^{1/2}$ times a constant. See Bass and Burdzy \cite{BassBurdzy},
Davis \cite{Davis1983}, Davis and Perkins \cite{DavisPerkins},
Greenwood and Perkins \cite{GreenwoodPerkins},
Kahane \cite{Kahane1974,Kahane1976,Kahane1983},
and Perkins \cite{Perkins1983}.
And item $(\mathscr{b})$ parallels
the assertion that the slowest points for the increments of Brownian motion
move as $\Omega(t^{1/2})$  at time $t\approx 0$; see Dvoretzky
\cite{Dvoretzky}.\footnote{As is customary, we write
	$f(t)=\Omega(g(t))$ for $t\approx0$ when we mean that
	$f(t)\gtrsim g(t)$ for all $t\in[0\,,t_0]$ for a sufficiently small $t_0>0$.}

Part (3) of Theorem \ref{th:BC:u} too has counterparts 
in the earlier literature  on the increments of
Brownian motion; see Davis and Perkins \cite{DavisPerkins}
and Perkins \cite{Perkins1983}. 
Moreover, just as in the theory of Brownian slow points,  there is a sharp
cutoff point, here at $\theta=\theta_c$; see \eqref{theta_c}. 

Among other things, Theorem \ref{th:BC:u}(3)
implies that the collection of points of slow growth -- that  is
$\mathfrak{S} := \cup_{\theta\ge\theta_c}\mathfrak{S}(\theta)$ --
has a multifractal structure in the sense that we can decompose $\mathfrak{S}$
as a countable union of disjoint
random sets $\mathfrak{S}_0,\mathfrak{S}_1,\mathfrak{S}_2,\ldots$ that satisfy
\begin{equation}\label{dimdec}
	\dimh \mathfrak{S}_n < \dimh \mathfrak{S}_{n+1}
	\quad\text{a.s.\ for every $n\in\Z_+$.}
\end{equation}
For example, one can see this by setting
$\mathfrak{S}_n := \mathfrak{S}(2^n\theta_c)\setminus \mathfrak{S}(2^{n-1}\theta_c)$
for every $n\in\N.$
Theorem \ref{th:BC:u}(3), and the strict monotocity of $\lambda$
[see Theorem \ref{th:BC:u}(1)], together imply that $\dimh \mathfrak{S}_n
=\dimh \mathfrak{S}(2^n\theta_c) = 1 - 2\lambda(2^n\theta_c)$ for
every $n\in\Z_+$,
whence follows our claim \eqref{dimdec}.

\begin{conjecture}
	In analogy with the theory of slow points for Brownian motion,
	we conjecture
	that $\fS(\theta_c)=\varnothing$ a.s.; compare
	with Bass and Burdzy \cite{BassBurdzy}. The best we currently know
	is that $\dimh\mathfrak{S}(\theta_c)=0$; see Theorem \ref{th:BC:u}(3).
\end{conjecture}

We have no idea how to evaluate the function $\lambda$ either exactly,
or numerically. Moreover, we believe that exact values of $\lambda$ are likely to
remain unknown, though  the last section of our recent paper
\cite{KL} can serve as
a beginning foray into a future possibility of Monte-Carlo evaluations of $\lambda$.
Despite these comments, it is possible to extract the following asymptotic information about the function
$\lambda$.

\begin{corollary}\label{cor:lambda}
	The function $\lambda$ satisfies the following:
	\begin{align*}
		-\infty< \liminf_{\theta\to\infty} \theta^{-2}\log\lambda(\theta)
			&\le \limsup_{\theta\to\infty} \theta^{-2}\log\lambda(\theta)<0;
			\quad\text{and}\\
		0<\liminf_{\theta\to0^+} \theta^{4}\lambda(\theta)
			&\le\limsup_{\theta\to0^+}\theta^{4}\lambda(\theta)<\infty.
	\end{align*}
\end{corollary}

We have so far presented many of the analogies between the study of the points
of slow growth for \eqref{u} and the theory of Brownian slow points. There
also are significant differences, mainly in the respective proofs, caused by
the inherently infinite-dimensional nature of the problem studied here:
The Brownian theory uses finite-dimensional Markov process theory in an 
essential way, whereas $t\mapsto u(t\,,x)$ is not Markovian for any fixed finite collection of
spatial points $x\in\R$.

We conclude the Introduction with the following, which is motivated by 
its counterpart in Theorem 1.1, and especially Eq.~(1.4), of 
Davis and Perkins \cite{DavisPerkins}.

\begin{conjecture}
	We believe that for every nonrandom compact $K\subset\R$,
		\[
			\P\{\fS(\theta)\cap K\neq\varnothing\}=
			\begin{cases}
				1 & \text{if $\lambda(\theta) < \frac12\dimh K$},\\
				0 & \text{if $\lambda(\theta) > \frac12\dimh K$}.
			\end{cases}
		\]
\end{conjecture}

Throughout this paper, and for reasons that have been explained already
in Remark \ref{rem:meas},
we assume that the 
underlying probability space $(\Omega\,,\mathscr{F}\,,\P)$
is complete. On a few occasions, we will use the notation,
$\log_+(a)=\log(a+\e)$ for all $a\ge0$.

\section{Linearization}

Recall that the solution to \eqref{u} is a predictable random field
$u=\{u(t\,,x)\}_{t\ge0,x\in\R}$,
in the sense of Walsh \cite{Walsh},
that satisfies $\sup_{t\in(0,T)}\sup_{x\in\R}\E(|u(t\,,x)|^2)<\infty$ for all $T>0$ and 
solves
\begin{equation}\label{mild:u}
	u(t\,,x) = 1 + \int_{(0,t)\times\R} G_{t-s}(y-x)\sigma(u(s\,,y))\,W(\d s\,\d y)
	\qquad\forall t>0,\ x\in\R,
\end{equation}
almost surely, where 
\begin{equation}\label{G}
	G_s(y) = (4\pi s)^{-\frac{1}{2}}\exp\left( - \frac{y^2}{4s}\right)\qquad\forall s > 0, y \in \R.
\end{equation}
Ordinarily, one would need to be careful about the order of
the quantifiers in \eqref{mild:u} 
[$\forall t$ and $\forall x$ and the ``almost sure'' portion of the assertion].
Such distinctions do not exist in the present context since, as was pointed out earlier
in the Introduction, $u$ is continuous.

Theorem \ref{th:BC:u} is a statement about the behavior of the infinite-dimensional
process $t\mapsto u(t)$ when $t\approx0$. It has been observed by
Khoshnevisan, Swanson, Xiao, and Zhang \cite{KSXZ}
and Hairer and Pardoux \cite{HairerPardoux}
that, because of \eqref{mild:u} and to leading order,
the small-$t$ behavior of $u$ is essentially the same as the small-$t$ behavior of
$\sigma(1)H$, where $H=\{H(t\,,x)\}_{t\ge0,x\in\R}$ solves \eqref{u} starting from zero
and with $\sigma$ replaced identically by one. Stated somewhat more carefully (but
still not quite rigorously),
\begin{equation}\label{loclin}
	u(t\,,x) = 1 + \sigma(1)H(t\,,x) + \mathscr{o}(t^{1/4})
	\quad\text{as $t\to0^+$},
\end{equation}
where
\begin{equation}\label{H}\left[\begin{aligned}
	&\partial_t H(t\,,x) = \partial^2_x H(t\,,x) + \dot{W}(t\,,x) &\qquad
		\text{for $(t\,,x)\in(0\,,\infty)\times\R$},\\
	&H(0) =  0&\text{on $\R$}.
\end{aligned}\right.\end{equation}
Just as in \eqref{mild:u}, we have
\begin{align}\label{mild:H}
	H(t\,, x) = \int_{(0, t) \times \R} G_{t-r}(z-x)\, W(\d r\, \d z)
	\qquad\forall t>0,\ x\in\R,
\end{align}
and $H(0\,,x)=0$ $[x\in\R$]. In particular, $H=\{H(t\,,x)\}_{t\ge0,x\in\R}$ is a 
centered Gaussian process, and the bulk of the proof of Theorem \ref{th:BC:u}
will involve a detailed analysis of the process $H$. 

The main purpose of this section is to
verify a suitably strong version of \eqref{loclin}.
The subsequent analysis of $H$, and its relation to Theorem \ref{th:BC:u},
will be developed in later sections.
First, let us record some of the salient features of the process $H$.

\begin{lemma}\label{lem:H}
	$H$ is a centered Gaussian process such that: 
	\begin{compactenum}
		\item $\Var H(t\,,x)
			=\sqrt{t/(2\pi)}$ for all $t\ge0$ and $x\in\R$.
		\item (H\"older continuity) $H\in C^{a/2,a}(\R_+\times\R)$ a.s.\ for every
			$a\in(0\,,\frac12)$.
		\item (The Markov property) The $C(\R)$-valued process $\{H(t)\}_{t\ge0}$
			is Markov.
		\item (Stationarity) $x\mapsto H(\cdot\,,x)$ is a $C(\R_+)$-valued stationary process.
		\item (Scaling) The law of
			$\{ c^{-1/4}H(ct\,,c^{1/2} x)\}_{t\ge0,x\in\R}$ does not depend on $c>0$.
	\end{compactenum}
\end{lemma}

Lemma \ref{lem:H} is essentially entirely a ready corollary of the following
formula, valid for every $s,t>0$ and $x,y\in\R$:
\begin{equation}\label{Cov:H}
	\Cov[H(t\,,x)\,,H(s\,,y)] 
	=\int_0^{s\wedge t}\exp\left( -
	\frac{|x-y|^2}{4(t+s-2r)}\right)
	\frac{\d r}{\sqrt{4\pi(t+s-2r)}}.
\end{equation}
The computation is straightforward. Therefore, we omit the details.
Let us also observe the following technical consequence of
\eqref{Cov:H}: 
\begin{equation}\label{H-H}\begin{split}
	\sup_{t>0}\E\left(|H(t\,,x)-H(t\,,y)|^2\right) &\asymp |x-y|,\\
	\sup_{t>0}\left( | H(t+\varepsilon\,,x) - H(t\,,x)|^2\right)
		&\asymp \sqrt{\varepsilon},
\end{split}\end{equation}
valid uniformly for all $\varepsilon>0$
and $x,y\in\R$; see \cite[\S3]{CBMS}.
In order to streamline our discussion, let us define 
\begin{align}\label{E(t,x)}
	\mathscr{E}(t\,,x) = u(t\,,x) - 1 - \sigma(1)H(t\,,x), \quad \forall t \ge 0, x \in \R,
\end{align}
where $u$ and $H$ respectively solve \eqref{u}
and \eqref{H}, and recall from \eqref{loclin} that
we aim to prove that $\mathscr{E}(t\,,x) = \mathscr{o}(t^{1/4})$
in a strong-enough sense when $t\approx0$. It turns out to be enough to 
carry out this plan in the seemingly special case that $\sigma$ is
in addition a bounded function. In that case,
we first prove a quantitative variation of
the improved pointwise statement that 
$\mathscr{E}(t\,,x) = \mathcal{O}(t^{1/2})$ in $L^k(\Omega)$
for every $k\ge1$.
Denote the $L^k(\Omega)$-norm of a random variable $X$ by $\|X\|_k = \E(|X|^k)^{1/k}$.

\begin{lemma}\label{lem:L_k:E}
	If, in addition, $\sigma$ is bounded, then
	$\sup_{x\in\R}\|\mathscr{E}(t\,,x)\|_k \lesssim k\sqrt{t}$
	uniformly for all $k \in [2\,,\infty)$ and $t \in (0\,,1]$.
	In particular, there exists $c>0$ such that
	\[
		\adjustlimits\sup_{t\in(0,1]}
		\sup_{x\in\R}\E\exp\left( c\, t^{-1/2}|\mathscr{E}(t\,,x)|\right)<\infty.
	\]
\end{lemma}

\begin{proof}
	The second, displayed, assertion of the lemma follows from the first 
	and the Taylor expansion of the exponential function. Therefore, it suffices
	to prove the asserted $L^k$-bound.
	
	We can combine \eqref{mild:u} and \eqref{mild:H}
	with a suitable formulation of the
	Burkholder-Davis-Gundy inequality for stochastic convolutions
	\cite[Proposition 5.2]{CBMS} in order to see that
	\begin{align*}
		\|\mathscr{E}(t\,,x)\|_k^2 &\le 4k \int_0^t \d s
			\int_{-\infty}^\infty\d y\ [G_{t-s}(x-y)]^2
			\|\sigma(u(s\,,y))-\sigma(1)\|_k^2\\
		&\le 4k\lip(\sigma)^2 \int_0^t\d s \int_{-\infty}^\infty\d y\ 
			[G_{t-s}(x-y)]^2
			\|u(s\,,y)-u(0\,,y)\|_k^2,
	\end{align*}
	uniformly for all $k \in [2\,,\infty)$, $t \in (0\,,1]$, and $x \in \R$.
	We now combine \eqref{mild:u} with \eqref{mild:H}
	and a suitable form of the Burkholder-Davis-Gundy inequality
	\cite[Proposition 5.2]{CBMS}
	in order to deduce from \eqref{H-H} that
	\begin{align}
		&\|u(t\,,x)-u(s\,,y)\|_k^2\notag \\
		& \le 4k \int_0^1 \d r \int_{-\infty}^\infty \d z \ 
			[G_{t-r}(x-z) \1_{(0,t)}(r) - G_{s-r}(y-z) \1_{(0,s)}(r)]^2 
			\|\sigma(u(r\,,z))\|_k^2\notag \\
		&\le 4k \sup_{z \in \R} |\sigma(z)|^2  
			\| H(t\,,x)-H(s\,,y)\|_2^2 \lesssim k \left(\sqrt{|t-s|} + |x-y|\right),
		\label{u-u}
	\end{align}
	uniformly for all $k \in [2\,,\infty)$, $t,s \in (0\,,1]$ 
	and $x, y \in\R$. Therefore, the preceding discussion, and the semigroup property
	of the heat kernel, together yield
	\[
		\|\mathscr{E}(t\,,x)\|_k^2 
		\lesssim k^2 \int_0^t\d s \int_{-\infty}^\infty
		\d y\ [G_{t-s}(x-y)]^2\sqrt{s}
		= \frac{k^2}{\sqrt{4\pi}}
		\int_0^t\sqrt{\frac{s}{t-s}}\,\d s
		\propto k^2 t,
	\]
	where the implied constant is independent of $(k\,,t\,,x)$.
	The lemma follows.
\end{proof}

Next we present a Gaussian upper bound
for the modulus of continuity of the two-parameter random field
$\mathscr{E}$.

\begin{lemma}\label{lem:E:exp:E}
	Suppose, in addition, that $\sigma$ is bounded, and choose
	and fix $N\in\N$. Then, there exists $\gamma_0>0$ such that
	\begin{align*}
		A= \sup_{\varepsilon\in(0,1/\e)}\E\left[ \exp\left( \gamma_0 \adjustlimits
		\sup_{\substack{ 0<s, t \le 1\\|t-s| \le\varepsilon^2}}
		\sup_{\substack{ -N\le x, y\le N\\|x-y|\le\varepsilon}}
		\frac{|\mathscr{E}(t\,,x) - \mathscr{E}(s\,,y)|^2}{%
		\varepsilon\log(1/\varepsilon)}
		\right) \right] < \infty.
	\end{align*}
	In particular, uniformly for all $\varepsilon\in(0\,,1/\e)$,
	\[
		\left\| \sup_{\substack{ 0<s, t \le 1\\|t-s| \le\varepsilon^{1/2}}}
		\sup_{\substack{ -N\le x, y\le N\\|x-y|\le\varepsilon}}
		|\mathscr{E}(t\,,x) - \mathscr{E}(s\,,y)|\right\|_k
		\le A \sqrt{k\varepsilon\log(1/\varepsilon)}.
	\]
\end{lemma}

\begin{proof}
	Recall \eqref{E(t,x)}.
	For every $t,s \in (0\,,1]$ and $x, y \in [-N\,,N]$, we may write
	\[
		\mathscr{E}(t\,,x)-\mathscr{E}(s\,,y) = u(t\,,x)-u(s\,,y) - \sigma(1)(H(t\,,x)-H(s\,,y)).
	\]
	Let $\Delta((t\,,x)\,,(s\,,y)) = |t-s|^{1/4} + |x-y|^{1/2}$ for all $s,t>0$
	and $x,y\in\R$.
	A crude application of \eqref{u-u}, once for $u$ and once for $H$, yields
	\begin{align}\begin{split}\label{E-E}
		&\|\mathscr{E}(t\,,x)-\mathscr{E}(s\,,y)\|_k \\
		&\lesssim 
			\|u(t\,,x)-u(s\,,y)\|_k + \|H(t\,,x)-H(s\,,y)\|_k
		\lesssim\sqrt{k}\,\Delta((t\,,x),(s\,,y)),
	\end{split}\end{align}
	uniformly for all $k \in [2\,,\infty)$, $t,s \in (0\,,1]$ and $x, y \in\R$. 
	Therefore, a standard metric entropy argument yields the sub-Gaussian 
	bound,
	\begin{align}\label{E:sup:E-E}
		\E\left[ \exp\left( \gamma_0 
		\sup_{\substack{ 0<s\neq t \le 1\\-N\le x\neq y\le N}}
		\frac{|\mathscr{E}(t\,,x) - \mathscr{E}(s\,,y)|^2}{%
		\left| \Delta ((t\,,x),(s\,,y))\right|^2
		\log_+\left(\dfrac{1}{\Delta((t\,,x),(s\,,y))} \right) } 
		\right) \right] < \infty,
	\end{align}
	valid for a suitably small choice of $\gamma_0>0$.
	For instance, \eqref{E:sup:E-E} can be obtained from \eqref{E-E} and 
	the Garsia-Rodemich-Rumsey lemma \cite{GRR} in a suitable form 
	(e.g., \cite[Proposition A.1]{DKN07} with $p(r)=r$ and 
	$\Psi(r) = \exp(c r^2)-1$ for $c>0$ small).
	The first assertion of the lemma follows readily from \eqref{E:sup:E-E}.	
	The second statement of
	the lemma follows from the first and the fact that if $X$ is a random
	variable such that $B=\E\exp(X^2)<\infty$, then
	\[\textstyle
		\| X\|_{2m}^{2m}/m^m \le
		\| X\|_{2m}^{2m}/m! \le \E\sum_{n=0}^\infty
		|X|^{2n}/n! = B,
	\]
	for every integer $m\ge1$.
\end{proof}

Lemma \ref{lem:E:exp:E} can now be used in order to improve
Lemma \ref{lem:L_k:E}, and yield the following quantitative uniform
improvement to \eqref{loclin}:

\begin{proposition}\label{pr:L_k:supE}
	Suppose in addition that $\sigma$ is bounded,
	and choose and fix two real numbers $a<b$. Then,
	uniformly for all $k \in [1\,,\infty)$ and $t \in (0\,,1/\e]$,
	\[
		\left\| \adjustlimits\sup_{s\in(0,t]}
		\sup_{x \in [a,b]} |\mathscr{E}(s\,,x)| \right\|_k \lesssim k \sqrt{t}|\log t|.
	\]
\end{proposition}

\begin{proof}
	We prove the proposition by appealing to an
	interpolation argument. Thanks to Jensen's inequality,
	it suffices to consider $k\ge 2$.
	
	Choose and fix two real numbers
	$a<b$. Also, let us write
	\[
		L_t = \sqrt{t}|\log t|\qquad\forall t\in(0\,,1/\e],
	\]
	in order to simplify the exposition in a few spots.
	Because both $u$ and $H$ are continuous, 
	$\sup_{s\in(0,t]}\sup_{x\in[a,b]}|\mathscr{E}(s\,,x)|$
	and other such objects discussed here are all random variables. Therefore,
	we may proceed without concerns for measurability issues.
	
	Define $J_n = \cup_{j=1}^{n^2}\{ j/n^2\}$ and
	$K_n = \cup_{j=0}^n \{ a + j(b-a)n^{-1}\}$ for all $n\in\N$.
	We can then write
	\[
		\P\left\{ \adjustlimits\sup_{s\in(0,t]}
		\sup_{x\in[a,b]} |\mathscr{E}(s\,,x)| \ge z\right\} 
		\le T_1 + T_2,
	\]
	where
	\begin{align*}
		T_1 &= \P\left\{ \adjustlimits\max_{r\in J_n}
			\max_{x\in K_n} |\mathscr{E}(r t \,,x)| \ge
			z/2\right\},\\
		T_2 &= \P\left\{\sup_{\substack{s,t\in (0,1/\e]:\\
			|t-s|\le 1/(\e n^2)}}
			\sup_{\substack{x,y\in [a,b]:\\
			|x-y|\le (b-a)/n}} |\mathscr{E}(s\,,x)-\mathscr{E}(t\,,y)| \ge
			z/2 \right\}.
	\end{align*}
	Since the respective cardinalities of $J_n$ and $K_n$ are $n^2$ and $n+1$, 
	Lemma \ref{lem:L_k:E} and Boole's inequality together yield
	a number $c_0>0$ such that
	\[
		T_1 \le n^2(n+1)\adjustlimits\sup_{s\in(0,t]}\sup_{x\in\R}
		\P\left\{ |\mathscr{E}(s\,,x)| 
		\ge z/2\right\} \lesssim n^3 \e^{-c_0z/\sqrt{t}},
	\]
	uniformly in $(t\,,z\,,n)\in(0\,,1/\e]\times(0\,,\infty)\times\N$. 
	Let $c_1 = \max\{ 1/\sqrt{\e}\,,b-a \}$.
	Then Lemma \ref{lem:E:exp:E} and Chebyshev's inequality together yield
	a number $c_2>0$ such that
	\begin{align*}
		T_2 &\le
			\P\left\{\sup_{\substack{s,t\in (0,1/\e]:\\|t-s|\le1/(\e n^2)}}
			\sup_{\substack{x,y\in [a,b]:\\
			|x-y|\le (b-a)/n}} \frac{|\mathscr{E}(s\,,x)-\mathscr{E}(t\,,y)|^2}{
			\left(\dfrac{c_1}{n}\right)\log_+\left(\dfrac{n}{c_1}\right)} \ge
			\frac{z^2}{4\left(\dfrac{c_1}{n}\right)
			\log_+\left(\dfrac{n}{c_1}\right)} \right\}\\
		&\lesssim\exp\left( - \frac{c_2nz^2}{\log_+ n}\right),
	\end{align*}
	uniformly in $(t\,,z\,,n)\in(0\,,1/\e]\times(0\,,\infty)\times\N$.
	Now, we may integrate
	by parts as follows: For every $A>0$,
	\begin{align*}
		\Upsilon_{t,k} &:= \E\left(\adjustlimits\sup_{s\in(0,t]}
			\sup_{x\in[a,b]} |\mathscr{E}(s\,,x)|^k\right) =
			k\int_0^\infty z^{k-1}\P\left\{ \sup_{s\in (0,t]}
			\sup_{x\in [a,b]} |\mathscr{E}(s\,,x)| \ge z\right\} \d z\\
		&\le (AL_t)^k + 
			k\int_{AL_t}^\infty z^{k-1} T_1 \,\d z
			+ k\int_{AL_t}^\infty z^{k-1} T_2\,\d z,
	\end{align*}
	pointwise. Thus, we plug in the above estimates for $T_1$ and $T_2$
	in order to see that
	\[
		\Upsilon_{t,k} \lesssim (AL_t)^k
		+  n^3k\int_{AL_t}^\infty z^{k-1}\e^{-c_0z/\sqrt{t}}\,\d z
		+ k\int_{AL_t}^\infty z^{k-1}
		\exp\left( - \frac{c_2nz^2}{\log_+ n}\right)\d z,
	\]
	uniformly in $(A\,,t\,,n\,,k)\in(0\,,\infty)\times(0\,,1/\e]\times\N\times\N$.
	Now,
	\begin{align*}
		&\int_{AL_t}^\infty z^{k-1}\e^{-c_0z/\sqrt{t}}\,\d z
			\le t^{k/2} \int_{A|\log t|}^\infty y^{k-1}\e^{-c_0y}\,\d y\\
		&\le t^{(k/2)+(Ac_0/2)} \int_{A|\log t|}^\infty y^{k-1}\e^{-c_0y/2}\,\d y
			\le (2/c_0)^k\Gamma(k)t^{(k/2)+(Ac_0/2)},
	\end{align*}
	uniformly in $(A\,,t\,,n\,,k)\in(0\,,\infty)\times
	(0\,,1/\e]\times\N\times\N$.
	Similarly,
	\[
		\int_{AL_t}^\infty z^{k-1}
		\exp\left( - \frac{c_2nz^2}{\log_+ n}\right)\d z
		= \left( \frac{\log_+ n}{n}\right)^{k/2}
		\int_{AL_t \sqrt{n/\log_+n}}^\infty y^{k-1} \e^{-c_2y^2}\,\d y,
	\]
	uniformly in $(A\,,t\,,n\,,k)\in(0\,,\infty)\times
	(0\,,1/\e]\times\N\times\N$. There exists $c_*>1$
	such that $L_t \sqrt{n/\log_+n} \ge 1$ for $n= n_t = \lfloor c_*/t\rfloor.$
	For this choice of $n=n_t$, we write
	\begin{align*}
		\int_{AL_t}^\infty z^{k-1}
			\exp\left( - \frac{c_2n_t z^2}{\log_+ n_t}\right)\d z
			&\le \left( \frac{\log_+ n_t}{n_t}\right)^{k/2}
			\int_0^\infty y^{k-1} \e^{-c_2y^2}\,\d y\\
		& \le (L_t/ \sqrt c_2)^{k}\Gamma(k/2),
	\end{align*}
	whence (for $n$ replaced by $n_t$),
	\[
		\Upsilon_{t,k} \lesssim (AL_t)^k
		+  k(2/c_0)^k\Gamma(k)t^{(k/2)+(Ac_0/2) -3}
		+ k(L_t/\sqrt c_2)^k\Gamma(k/2),
	\]
	uniformly in $(A\,,t\,,k)\in(0\,,\infty)\times
	(0\,,1/\e]\times\N$. We may choose
	$A=6/c_0$ in order to deduce the following:
	\[
		\Upsilon_{t,k}^{1/k}
		\lesssim \left( 1
		+  k\Gamma(k)
		+ k\Gamma(k/2)\right)^{1/k} L_t,
	\]
	uniformly in $(t\,,k)\in(0\,,1/\e]\times\N$. 
	An application of Stirling's formula
	implies the proposition.
\end{proof}

We now state and prove the analogue of Part (1) of
Theorem \ref{th:BC:u}, valid instead for the random
field $H$. Thanks to Lemma \ref{lem:H} the
law of the process $H(\cdot\,,x)$ does not depend on $x$.
Therefore, we write the following for $x=0$ to simplify the
notation a little.

\begin{proposition}\label{pr:BC:H}
	For every $\theta>0$,
	\[ 
		\P\left\{  | H(t\,,0)| \le \theta t^{1/4}\
		\forall t\in[a\,,b]\right\}
		= (a/b)^{\lambda(\theta)+\mathscr{o}(1)}\quad\text{as $a/b\to0$},
	\]
	where $\lambda:(0\,,\infty)\to(0\,,\infty)$ is strictly decreasing, convex, continuous,
	and satisfies $\lim_{\theta\to0^+}\lambda(\theta)=\infty$
	and $\lim_{\theta\to\infty}\lambda(\theta)=0$.
\end{proposition}

\begin{proof}
	Thanks to scaling (Lemma \ref{lem:H}),
	\[ 
		\P\left\{  | H(t\,,0)| \le \theta t^{1/4}\ \forall t\in[a\,,b]\right\}
		=  \P\left\{  | H(t\,,0)| \le \theta t^{1/4}\
		\forall t\in[1\,,b/a]\right\}.
	\]
	Therefore, it suffices to prove that
	\begin{align}\label{bcp}
		\P\left\{  | H(t\,,0)| \le \theta t^{1/4}\
		\forall t\in[1\,,N]\right\}
		= N^{-\lambda(\theta)+\mathscr{o}(1)}\quad\text{as $N\to\infty$}.
	\end{align}
	In fact, Theorem 1.1 of our paper \cite{KL} states that if 
	$\{X(t)\}_{t \ge 0}$ is a centered Gaussian process with continuous 
	sample functions that is self-similar with index $\alpha>0$ and 
	satisfies strong local non-determinism, then there exists $\lambda$ 
	with the asserted properties such that 
	$\P\{ |X(t)| \le c t^\alpha \ \forall t \in [1\,,T] \} = 
	T^{-\lambda(c)+\mathscr{o}(1)}$ as $T \to \infty$.
	Hence, \eqref{bcp} can be obtained by taking $X(t) = H(t\,,0)$ 
	and applying that theorem, whose condition are verified in Section 3 of \cite{KL},
	using parameters $d=\nu=\beta=\gamma=1$ and $\alpha=1/4$.
\end{proof}

\section{Proof of Theorem \ref{th:BC:u}, Part (1)}

We begin our work by deriving Theorem \ref{th:BC:u}(1)
from Proposition \ref{pr:BC:H}. It is possible to go into 
our estimates more deeply and improve the condition of
Theorem \ref{th:BC:u} on $f$, in case there is need. 
But we will not do that here. 

Our proof of Theorem \ref{th:BC:u} is carried out in a natural way
in three steps. Throughout, we choose and fix a number 
$\theta>0$, as in the statement of the theorem.\\

\noindent\emph{Step 1.} In this first step we prove Theorem \ref{th:BC:u}(1)
under the additional hypothesis that $\sigma$ is bounded. Let us therefore
assume that. 

Recall the linearization error $\mathscr{E}$ defined in \eqref{E(t,x)}.
According to Proposition \ref{pr:L_k:supE},
\[
	\left\| \sup_{s\in(\e^{-n},\e^{-n+1}]}
	\frac{|\mathscr{E}(s\,,0)|}{s^{1/4}}\right\|_k
	\le \e^{n/4} \left\| \sup_{s\in(\e^{-n},\e^{-n+1}]}
	|\mathscr{E}(s\,,0)|\right\|_k
	\lesssim n\exp(-n/4) k,
\]
uniformly for all integers $n\ge2$ and $k\ge1$. Consequently,
\[
	\left\| \sup_{s\in(0,\e^{-n+1}]}
	\frac{|\mathscr{E}(s\,,0)|}{s^{1/4}}\right\|_k
	\lesssim k\sum_{j=n}^\infty j \exp(-j/4)\lesssim kn \exp(-n/4),
\]
uniformly for all integers $n,k\ge1$. Thus, there exists $c_0,c_1>0$ such that
\[
	\left\| \sup_{s\in(0,\varepsilon]}
	\frac{|\mathscr{E}(s\,,0)|}{s^{1/4}}\right\|_k
	\le c_0 k n\exp(-n/4) \le c_1 k\varepsilon^{1/4}|\log\varepsilon|,
\]
uniformly for all  $n,k\in\N$ and $\varepsilon\in(\exp(-n/4)\,, \exp(-(n-1)/4)]$.
Therefore, if we choose $c\in(0\,,1/c_1)$, then
\begin{align*}
	\sup_{\varepsilon\in(0,1/\e]}
	\E\exp\left( \frac{c}{\varepsilon^{1/4}|\log\varepsilon|}\sup_{s\in(0,\varepsilon]}
	\frac{|\mathscr{E}(s\,,0)|}{s^{1/4}}\right)
	\lesssim 1 + \sum_{k=1}^\infty \frac{k^k(cc_1)^k}{k!}<\infty.
\end{align*}
In particular, Chebyshev's inequality yields
\[
	\P\left\{ \sup_{s\in(0,\varepsilon]}
	\frac{|\mathscr{E}(s\,,0)|}{s^{1/4}} \ge \delta |\sigma(1)|\right\}
	\lesssim \exp\left( - \frac{c\delta|\sigma(1)|}{\varepsilon^{1/4}|\log\varepsilon|}\right),
\]
uniformly for all $\varepsilon\in(0\,,1/\e]$ and $\delta>0$.
Therefore, for all $\delta>0$,
\begin{align*}
	&\P\left\{ \sup_{t \in [\varepsilon f(\varepsilon),\varepsilon]} 
		\frac{|u(t\,,0)-1|}{t^{1/4}} \le |\sigma(1)|\theta \right\}\\
	&\hskip.5in\lesssim \P\left\{ \sup_{t \in [\varepsilon f(\varepsilon),\varepsilon]}
		\frac{|H(t\,,0)|}{t^{1/4}} \le \theta + \delta \right\}
		+ \exp\left( - \frac{c\delta|\sigma(1)|}{\varepsilon^{1/4}|\log\varepsilon|}\right)\\
	&\hskip.5in= f(\varepsilon)^{\lambda(\theta+\delta)+\mathscr{o}(1)}
		+ \exp\left( - \frac{c\delta|\sigma(1)|}{\varepsilon^{1/4}|\log\varepsilon|}\right)
		= f(\varepsilon)^{\lambda(\theta+\delta)+\mathscr{o}(1)}
		\quad\text{as $\varepsilon\to0^+$},
\end{align*}
where the first identity in the last line holds thanks to scaling (Lemma \ref{lem:H})
and Proposition \ref{pr:BC:H}, and the second thanks to the condition 
of Theorem \ref{th:BC:u} on the function $f$.

Similarly, if $0<\delta<\theta$, then
\begin{align*}
	&\P\left\{ \sup_{t \in [\varepsilon f(\varepsilon),\varepsilon]}
		\frac{|u(t\,,0)-1|}{t^{1/4}} \le |\sigma(1)|\theta \right\} \\
	&\ge \P\left\{ \sup_{t \in [\varepsilon f(\varepsilon),\varepsilon]}
		\frac{|H(t\,,0)|}{t^{1/4}} \le \theta - \delta \right\} - 
		\P\left\{ \sup_{t \in [\varepsilon f(\varepsilon),\varepsilon]}
		\frac{|\mathscr{E}(t\,,0)|}{t^{1/4}} \ge \delta |\sigma(1)| \right\}\\
	&\gtrsim f(\varepsilon)^{\lambda(\theta-\delta)+\mathscr{o}(1)} - 
		\exp\left( - \frac{c\delta|\sigma(1)|}{\varepsilon^{1/4}|\log\varepsilon|}\right)
		= f(\varepsilon)^{\lambda(\theta-\delta) + \mathscr{o}(1)},
\end{align*}
as $\varepsilon\to 0^+$. Combine the preceding efforts in order to see that,
as $\varepsilon\to0^+$,
\[
	\lambda(\theta+\delta) +\mathscr{o}(1) \le 
	\frac{1}{\log f(\varepsilon)}\log\P\left\{ \sup_{t \in [\varepsilon
	f(\varepsilon),\varepsilon]} 
	\frac{|u(t\,,0)-1|}{t^{1/4}} \le |\sigma(1)|\theta \right\}\le 
	\lambda(\theta-\delta) +\mathscr{o}(1),
\]
for every fixed $0<\delta<\theta$.
Because $\lambda$ is continuous (Proposition \ref{pr:BC:H}), we may 
let $\delta \to 0^+$ to conclude the proof of Theorem \ref{th:BC:u}(1) under the additional
hypothesis that $\sigma$ is bounded.\\

\noindent\emph{Step 2.} In this step we show that the solution to \eqref{u}
locally in time behaves as the solution to the same SPDE but with a truncated
$\sigma$. Define
\[
	\bar\sigma(x) = \begin{cases}
		\sigma(2)&\text{if $x>2$},\\
		\sigma(x)&\text{if $|x|\le2$},\\
		\sigma(-2)&\text{if $x<-2$},
	\end{cases}
\]
and let $\bar{u}$ solve the initial-value problem \eqref{u} using
the same noise $\dot{W}$ as before, but
with $\sigma$ replaced by $\bar\sigma$. Thanks to \eqref{mild:u}
and a suitable version of the Burkholder-Davis-Gundy inequality \cite[Proposition 5.2]{CBMS},
for all $k\in[2\,,\infty)$, $t>0$, and $x\in\R$,
\begin{align}
	\| u(t\,,x) - \bar{u}(t\,,x) \|_k^2
	&\le 4k\int_0^t\d s\int_{-\infty}^\infty\d y\
	|G_{t-s}(y-x)|^2 \left\| \sigma(u(s\,,y)) - \bar\sigma(\bar{u}(s\,,y))
	\right\|_k^2 \notag\\
	&\le T_1 + T_2,
	\label{u-bar(u)}
\end{align}
where
\begin{align*}
	T_1 &= 4k\lip(\sigma)^2\int_0^t\d s\int_{-\infty}^\infty\d y\
		|G_{t-s}(y-x)|^2 \left\| u(s\,,y) -\bar{u}(s\,,y)
		\right\|_k^2,\\
	T_2 &= 4k\int_0^t\d s\int_{-\infty}^\infty\d y\
		|G_{t-s}(y-x)|^2 \left\| \left\{ |\sigma(u(s\,,y))| + |\bar\sigma(\bar{u}(s\,,y))|\right\}
		\1_{\{\|\bar{u}(s\,,y)\|>2\}}\right\|_k^2.
\end{align*}
Though a great deal more is known,
we will use only the basic facts that
\[
	L_k = \adjustlimits\sup_{s\in[0,1]}\sup_{y\in\R} 
	\E\left( |\sigma(u(s\,,y))|^k+|\bar{\sigma}(\bar{u}(s\,,y))|^k\right)
\]
and
\[
	\bar{L}_k = \adjustlimits\sup_{t\in[0,1]}\sup_{x\in\R}
	t^{-k/4} \E\left( | \bar{u}(t\,,x) - 1|^k\right)
\]
are finite for every $k\ge2$; see Khoshnevisan \cite[Lemma 5.4 and Theorem 5.5]{CBMS}. Define
\[
	M(r) = \sup_{t \in [0,r]} 
	\sup_{x\in\R} \|u(t\,,x) - \bar{u}(t\,,x) \|_k^2
	\qquad\forall r\in[0\,,1],
\]
and observe that, for every $r \in (0\,,1]$ and $t\in(0\,,r]$,
\[
	T_1 \lesssim M(t)\int_0^t\d s\int_{-\infty}^\infty\d y\ |G_s(y)|^2 	
	= M(t)\int_0^t\frac{\d s}{\sqrt{4\pi s}}\propto M(t)\sqrt{t} \le M(r) \sqrt{r},
\]
where the implied constant does not depend on $r \in (0\,,1]$ or $t\in(0\,,r]$. Also,
by the Cauchy-Schwarz inequality,
\begin{align*}
	T_2 &\lesssim 4k\int_0^t\d s\int_{-\infty}^\infty\d y\
		|G_{t-s}(y-x)|^2 \left\|  |\sigma(u(s\,,y))| + |\bar\sigma(\bar{u}(s\,,y))| \right\|_{2k}^2
		\P\left\{\|\bar{u}(s\,,y)\|>2\right\}^{1/k}\\
	&\le 4kL_{2k}^{1/k}\int_0^t\d s\int_{-\infty}^\infty\d y\
		|G_{t-s}(y-x)|^2  	\P\left\{\|\bar{u}(s\,,y)-1\|>1\right\}^{1/k}\\
	&\le 4k (\bar{L}_kL_{2k})^{1/k}\int_0^t\d s\int_{-\infty}^\infty\d y\
		|G_{t-s}(y-x)|^2 s^{1/4}\\
		& =\frac{2k (\bar{L}_kL_{2k})^{1/k}}{\sqrt\pi}
		\int_0^t \frac{s^{1/4}}{\sqrt{t-s}}\,\d s
		\propto t^{3/4} \le r^{3/4},
\end{align*}
where the implied constant does not depend on $r \in (0\,,1]$ or $t\in(0\,,r]$.
It follows from the above estimates for $T_1$ and $T_2$, and
from \eqref{u-bar(u)}, that $\|u(t\,,x)-\bar u(t\,,x)\|_k^2 \lesssim M(r) \sqrt{r} + r^{3/4}$ uniformly for all $r \in (0\,,1]$, $t \in (0\,,r]$ and $x \in \R$.
Then, taking $\sup$ over all $t \in (0\,,r]$ and $x \in \R$ yields
$M(r) \lesssim M(r) \sqrt{r} + r^{3/4}$ uniformly for all $r\in (0\,,1]$.
This implies that $M(r)\lesssim r^{3/4}$ uniformly for all small-enough
$r\in(0\,,1]$ and hence for all $r\in(0\,,1]$.
This shows in particular that
\begin{equation}\label{uu}
	C_k =\sup_{t\in(0,1]} t^{-3k/8}
	\|  u(t\,,0) - \bar{u}(t\,,0) \|_k^k <\infty
	\quad\forall k\ge2.
\end{equation}
We will need the following variation which has the supremum inside the expectation.

\begin{lemma}\label{lem:u-bar(u)}
	For every $k\ge 2$ and $\nu\in(0\,,3/8)$,
	\[
		\E\left( \sup_{t\in(0,1]}\frac{|u(t\,,0) - \bar{u}(t\,,0)|^k}{t^{k\nu}}\right)<\infty.
	\]
\end{lemma}

\begin{proof}
	Throughout, we choose and fix two numbers 
	$\nu\in (0\,,\frac38)$
	and $\delta \in (0\,,\frac38-\nu)$.
	For all $n \in \Z_+$ and $k \in [2/\delta\,,\infty)$, let $F_{n,k}$ 
	denote an equally-spaced mesh in $I(n) = [1/(n+1)\,,1/n]$,
	with mesh size $n^{-\delta k}$.
	Then, for every $z > 0$, $n \in \Z_+$ and $k \in [2/\delta\,,\infty)$, 
	we may use interpolation to write
	\begin{equation}\label{j123}
		\P\left\{ \sup_{t \in I(n)}
		\frac{|u(t\,,0)-\bar{u}(t\,,0)|}{t^\nu} \ge z \right\} \le J_1 + J_2 + J_3,
	\end{equation}
	where
	\begin{align*}
		&J_1 = \P\left\{ \max_{t\in F_{n,k}} 
			\frac{|u(t\,,0)-\bar{u}(t\,,0)|}{t^\nu} \ge \frac z3 \right\},\\
		&J_2 = \P\left\{ \sup_{t,s\in(0,1]: |t-s| \le n^{-\delta k}} 
			|\bar{u}(t\,,0)-\bar{u}(s\,,0)| \ge \frac{ z}{3(n+1)^\nu} \right\},\\
		&J_3 = \P\left\{ \sup_{t,s\in(0,1]: |t-s| \le n^{-\delta k}} 
			|u(t\,,0)-u(s\,,0)| \ge \frac{z}{3(n+1)^\nu} \right\}.
	\end{align*}
	We apply Boole's inequality, Chebyshev's inequality, and \eqref{uu}, in order to see that
	\begin{equation}\label{J1}
		J_1 \lesssim n^{-2} 3^k C_k n^{-k(\frac38-\nu-\delta)} z^{-k},
	\end{equation}
	uniformly for all $z>0$, $n \in \Z_+$ and $k \in [2/\delta\,,\infty)$.
	Since $\bar{\sigma}$ is bounded, we may use the 
	Burkholder-Davis-Gundy inequality \cite[Proposition 5.2]{CBMS}
	and \eqref{H-H} to deduce that
	\begin{align}\label{bar:u-u}
		\|\bar{u}(t\,,0)-\bar{u}(s\,,0)\|_k \lesssim \sqrt k |t-s|^{1/4},
	\end{align}
	uniformly for all $k \in [2/\delta\,,\infty)$ and $t,s \in (0\,,1]$.
	Hence, a standard metric entropy argument yields a number $c_0>0$ such that
	\[
		\E\left[ \exp\left( c_0 \sup_{0 < s < t \le 1} 
		\left|\frac{\bar{u}(t\,,0)-\bar{u}(s\,,0)}{|t-s|^{1/4}
		\sqrt{\log_+\frac{1}{|t-s|}}}\right|^2 \right) \right] < \infty,
	\]
	which follows, for instance, from \eqref{bar:u-u} and Proposition A.1 of 
	Dalang et al \cite{DKN07}.
	Thanks to Chebyshev's inequality and the preceding, there
	exists a number $L_1>0$ such that the following holds uniformly
	for all $z>0$, $n \in \Z_+$ and $k \in [2/\delta\,,\infty)$:
	\begin{equation}\label{J2}
		J_2 \le L_1 \exp\left( - \frac{n^{(\delta k/2) -2\nu} z^2}{L_1 k \log n} \right).
	\end{equation}
	Next, recall the following well-known estimate \cite[Theorem 5.5]{CBMS}: 
	There exists a number $L>0$ with the following property:
	For every $p\ge2$ there exists $L_*=L_*(p)>0$
	such that
	\begin{align}\label{u:lk}
		\|u(t\,,x)\|_p \le L_* \exp(L p^2 t) \quad
		\forall t>0,\ x\in\R.
	\end{align}
	Thanks to the Burkholder-Davis-Gundy inequality
	\cite[Proposition 5.2]{CBMS}, linear growth of $\sigma$, \eqref{u:lk}, 
	and \eqref{H-H}, there is a number $L_2>0$ such that uniformly
	for all $p \in [2/\delta\,,\infty)$ and $t,s\in (0\,,1]$,
	\begin{align}
		&\|u(t\,,0)-u(s\,,0)\|_p^2\notag\\
		 &\le 4 p \int_0^1 \d r \int_{-\infty}^\infty \d y \, 
		 	|G_{t-r}(y)\1_{(0,t)}(r)-G_{s-r}(y)\1_{(0,s)}(r)|^2 \|\sigma(u(r\,,y))\|_p^2\notag\\
		 &\lesssim \e^{L_2 p^2} |t-s|^{1/2}.\label{u-u:lk}
	\end{align}
	Choose and fix $L_3>L_2$ and define
	\[
		\Psi(x) = \sum_{p=1}^\infty \exp\{-L_3 (2p)^3\} x^{2p}
		\qquad\forall x\in\R.
	\]
	We may note that $\Psi$ is a strong Young function in the sense that:
	\begin{compactenum}
		\item $\Psi$ is even and convex on $\R$;
		\item It is strictly increasing on $\R_+$;
		\item $\Psi(0)=0$ and $\Psi(\infty)=\infty$; and
		\item $\Psi$ has a strictly increasing inverse $\Psi^{-1}$ on $\R_+$.
	\end{compactenum}
	Define
	\[
		\mathscr{C} = \int_0^1 \int_0^1 \Psi\left( \frac{|u(t\,,0)-u(s\,,0)|}{|t-s|^{1/4}} \right) \d s \, \d t.
	\]
	By Lemma 1.1 of Garsia, Rodemich, and Ramsey \cite{GarsiaRR},
	\begin{align*}
		|u(t\,,0)-u(s\,,0)| \le 2 \int_0^{|t-s|}
		\Psi^{-1}\left(\frac{4\mathscr{C}}{r^2}\right)
		\frac{\d r}{r^{3/4}} \qquad \forall t, s \in (0\,,1].
	\end{align*}
	By monotonicity, $\Psi^{-1}(y) \le \exp\{L_3 p^2\} y^{1/p}$ for every $y\ge0$ 
	and all even integers $p \ge 2$.
	Thus it follows that,
	for all even integers $p > 8$ there exists $L_4=L_4(p)>0$ such that
	\begin{align*}
		|u(t\,,0)-u(s\,,0)|
		\le L_4 \e^{L_3 p^2} \mathscr{C}^{1/p} |t-s|^{(p-8)/(4p)},
	\end{align*}
	uniformly for every $s,t\in[0\,,1]$.
	The preceding and \eqref{u-u:lk} together
	imply that for every even integer $p>8$ there exists $L_5 = L_5(p)>0$ such that
	\begin{align*}
		\E\left( \sup_{0<s< t \le 1} \frac{|u(t\,,0)-u(s\,,0)|^p}{|t-s|^{(p-8)/4}} \right)
			&\le L_4^p \e^{L_3 p^3}  \int_0^1 \int_0^1
			\sum_{m=1}^\infty \frac{\|u(t\,,0)-u(s\,,0)
			\|_{2m}^{2m}}{\e^{L_3(2m)^3}|t-s|^{m/2}} \d s \, \d t\\
		&\le L_5 \e^{L_3 p^3} \sum_{m=1}^\infty \e^{-(L_3-L_2) (2m)^3} < \infty.
	\end{align*}
	This together with Chebyshev's inequality yields
	\begin{equation}\label{J3}\begin{split}
		J_3 &\le \P\left\{ \sup_{0<s < t \le 1}
			\frac{|u(t\,,0)-u(s\,,0)|^p}{|t-s|^{(p-8)/4}} 
			\ge \frac{n^{\delta k(p-8)/4}  z^p}{(n+1)^{\nu p}3^p } \right\}
			\lesssim 3^p \e^{L_3 p^3} n^{\nu p-(\delta k (p-8)/4)} z^{-p}\\ 
		&= 3^p \e^{L_3 p^3} n^{-p(\delta k/4 - \nu - 2\delta k/p)} z^{-p},
	\end{split}\end{equation}
	uniformly for all $k \in [4/\delta\,,\infty)$ and all even numbers $p > 8$.
	Set 
	\[
		k_0 = 4(\nu+1)/\delta.
	\]
	If $k \in (k_0\,,\infty)$, then $(\delta k/4) - \nu > 1$.
	Thus, we may choose a large-enough even number $p = p(k\,,\delta) > \max\{ 8\,, k\}$ such that
	$p(\delta k/4 - \nu - 2\delta k/p) > p(1 - 2\delta)>1$
	We now return to \eqref{j123}, and apply \eqref{J1}, \eqref{J2},
	and \eqref{J3}, in order to see that for every $k\in(k_0\,,\infty)$ there exists
	$C>0$ such that
	\begin{align*}
		&\P\left\{ \sup_{t\in I(n)} \frac{|u(t\,,0)-\bar{u}(t\,,0)|}{t^\nu} \ge z \right\}\\
		&\lesssim  n^{-2-k(\frac38-\nu-\delta)} z^{-k}
		+  \exp\left( - \frac{n^{(\delta k/2) - 2\nu} z^2}{C \log n} \right)
		+ n^{-p(\delta k/4 - \nu - \frac{2\delta k}{p})} z^{-p},
	\end{align*}
	uniformly for all $z>0$ and $n \in \Z_+$.
	Sum the preceding over $n \in \Z_+$ in order to see that
	\[
		\sup_{z>0} z^k \P\left\{ \sup_{t \in (0,1]} 
		\frac{|u(t\,,0)-\bar{u}(t\,,0)|}{t^\nu} \ge z \right\} < \infty \qquad \forall k > k_0,
	\]
	whence follows the lemma.
\end{proof}

\noindent\emph{Step 3.}
We now conclude the proof of Part (1) of
Theorem \ref{th:BC:u}. As $\varepsilon\to0^+$,
\begin{align*}
	&\P\left\{ \sup_{t \in [\varepsilon f(\varepsilon),\varepsilon]} 
		\frac{|u(t\,,0)-1|}{t^{1/4}} \le |\sigma(1)|\theta \right\}\\
	&\lesssim \P\left\{ \sup_{t \in [\varepsilon f(\varepsilon),\varepsilon]}
		\frac{|\bar{u}(t\,,0)|}{t^{1/4}} \le  |\sigma(1)| (\theta + \delta) \right\}
		+\P\left\{ \sup_{t\in(\varepsilon f(\varepsilon),\varepsilon]} 
		\frac{|u(t\,,0) - \bar{u} (t\,,0)|}{t^{1/4}}\ge  |\sigma(1)|\delta \right\}\\
	&\le f(\varepsilon)^{\lambda(\theta+\delta)+\mathscr{o}(1)}
		+\P\left\{ \sup_{t\in(\varepsilon f(\varepsilon),\varepsilon]} 
		\frac{|u(t\,,0) - \bar{u} (t\,,0)|}{t^{1/2}}\ge
		\frac{ |\sigma(1)|\delta}{\varepsilon^{1/4}} \right\},
\end{align*}
thanks to Step 1, applicable since: (i)
$\bar\sigma$ is bounded as well as Lipschitz continuous; and
(ii) $\sigma(1)=\bar{\sigma}(1)\neq0$.
Therefore, Step 2 (specifically, Lemma \ref{lem:u-bar(u)}) ensures that for all $k\ge2$,
\begin{align*}
	\P\left\{ \sup_{t \in [\varepsilon f(\varepsilon),\varepsilon]} 
		\frac{|u(t\,,0)-1|}{t^{1/4}} \le |\sigma(1)|\theta \right\}
		&\lesssim f(\varepsilon)^{\lambda(\theta+\delta)+\mathscr{o}(1)}
		+\varepsilon^{k/4}\\
	&\le f(\varepsilon)^{\lambda(\theta+\delta)+\mathscr{o}(1)}
	\qquad\text{as $\varepsilon\to0^+$},
\end{align*}
thanks to the hypothesis of Theorem \ref{th:BC:u}
on the function $f$ and the fact that we can select $k$ to be as large as we want.
The very same argument can be recycled to show that
\[
	\P\left\{ \sup_{t \in [\varepsilon f(\varepsilon),\varepsilon]} 
	\frac{|u(t\,,0)-1|}{t^{1/4}} \le |\sigma(1)|\theta \right\}
	\ge f(\varepsilon)^{\lambda(\theta-\delta)+\mathscr{o}(1)}
	\qquad\text{as $\varepsilon\to0^+$},
\]
provided additionally that $\delta\in(0\,,\theta)$. This proves that
for all $\delta\in(0\,,\theta)$, and
as $\varepsilon\to0^+$,
\[
	\lambda(\theta+\delta)+\mathscr{o}(1)\le
	\frac{1}{\log f(\varepsilon)}
	\log\P\left\{ \sup_{t \in [\varepsilon f(\varepsilon),\varepsilon]} 
	\frac{|u(t\,,0)-1|}{t^{1/4}} \le |\sigma(1)|\theta \right\} \le
	\lambda(\theta-\delta)+\mathscr{o}(1).
\]
The continuity of $\lambda$ (see either Step 1 or Proposition \ref{pr:BC:H})
allows us to let $\delta$ tend to zero and
conclude the proof of part (1) of Theorem \ref{th:BC:u}.
\qed

\section{Proof of Corollary \ref{cor:lambda}}
In light of our earlier results in \cite{KL},
our proof of Corollary \ref{cor:lambda}
is brief. Indeed, we showed in Proposition \ref{pr:BC:H} that the function $\lambda$
is the so-called ``boundary-crossing exponent'' of the Gaussian
process $H$ in \eqref{H}. Lemma \ref{lem:H} and the material of
Section 3 of \cite{KL} together show that we may apply
Corollary 1.2 of \cite{KL} -- with $\delta=\alpha=\frac14$ -- in order
to deduce our Corollary \ref{cor:lambda}.\qed

\section{Localization}

One of the central ingredients of the proof of the remainder of
Theorem \ref{th:BC:u} is a localization property of $H$ from \eqref{H}.
Let us define for every $\alpha\in(0\,,1)$, $t>0$, and $x\in\R$,
the following ``localized'' version of the random field $H$.
\begin{equation}\label{H_alpha}
	H_\alpha(t\,,x) = \int_{(0,t)\times[x-t^{(1-\alpha)/2},\,x+t^{(1-\alpha)/2}]}
	G_{t-s}(y-x)\,W(\d s\,\d y).
\end{equation}
One can define also $H_\alpha(0)\equiv0$, essentially by continuity. 

The following is the main result of this section; it states that the effect of the noise
on the integral definition of $H$ in \eqref{H} is highly localized, primarily 
due to the rapid decay in the tails
of the heat kernel in \eqref{G}.

\begin{proposition}\label{pr:H-H_alpha:sup}
	Choose and fix some $\alpha\in(0\,,1)$ and $R>0$. Then, a.s.,
	\[
		\sup_{x\in[-R,R]}\left| H(t\,,x)-H_\alpha(t\,,x)\right| = \mathscr{o}(t^{1/4})
		\quad\text{as $t\to0^+$}.
	\]
\end{proposition}

This proposition holds because of its quantitative counterpart in Lemma 
\ref{lem:H-H_alpha:prob} below. That result, in turn, is derived in a series
of smaller steps, beginning with the following simple estimate.

\begin{lemma}\label{lem:H-H_alpha}
	For all $\alpha\in(0\,,1)$, $x\in\R$, and $t>0$,
	\[
		\E\left(|H(t\,,x) - H_\alpha(t\,,x)|^2\right)\le 
		2\sqrt{\frac{t}{\pi}} \exp\left( - \frac{1}{4t^\alpha}\right).
	\]
\end{lemma}

\begin{proof}
	By stationarity, the second moment in question does not depend on $x$;
	therefore, we consider only $x=0$. 
	Let $X$ denote a random variable with a standard normal distribution
	and recall that $\P\{|X|>y\}\le\exp(-y^2/2)$ for all $y>0$; see for example 
	\cite[Lemma A.3]{CBMS}.
	Since $G_s(y)\le  (4\pi s)^{-1/2}$ for all $s>0$ and $y\in\R$,
	it follows from the $L^2$-isometry for Wiener integrals that
	\begin{align*}
		& \E\left( |H(t\,,0)-H_\alpha(t\,,0)|^2\right)
			=  2\int_0^t\d s\int_{t^{(1-\alpha)/2}}^\infty\d y\
			|G_s(y)|^2\\
		& \le\int_0^t (\pi s)^{-1/2}\P\left\{ |X| > \sqrt{\frac{t^{1-\alpha}}{2s}} \right\}\d s
			\le\int_0^t (\pi s)^{-1/2}\exp\left(-\frac{t^{1-\alpha}}{4s}\right)\d s.
	\end{align*}
	This proves the lemma since $\exp\{-t^{1-\alpha}/(4s)\} \le\exp\{-t^{-\alpha}/4\}$ when $s\le t$.
\end{proof}

Next, we present a modulus of continuity estimate for the spatial variable of
the localized Gaussian random
field $H_\alpha$ in \eqref{H_alpha}.

\begin{lemma}\label{lem:H_alpha-H_alpha:x}
	Uniformly for all $x,y\in\R$, 
	\[
		\adjustlimits\sup_{t\in(0,1]}\sup_{\alpha\in(0,1)}
		\E\left(|H_\alpha(t\,,x) - H_\alpha(t\,,y)|^2\right) 
		\lesssim|x-y|\log_+\left(\frac{1}{|x-y|}\right).
	\]
\end{lemma}
\begin{proof}
	Without loss of generality, we may and will assume that
	$x>y$ throughout the proof.
	
	We can write
	\[
		H_\alpha(t\,,x) - H_\alpha(t\,,y) =I_1 + I_2,
	\]
	where
	\begin{align*}
		I_1 &= \int_{(0,t)\times[x-t^{(1-\alpha)/2},\,x+t^{(1-\alpha)/2}]}
			\left[G_{t-s}(z-x) - G_{t-s}(z-y)\right] W(\d s\,\d z),\\
		I_2 &=\int_{(0,t)\times[x-t^{(1-\alpha)/2},\,x+t^{(1-\alpha)/2}]
			\setminus [y-t^{(1-\alpha)/2},\,y+t^{(1-\alpha)/2}]}
			G_{t-s}(z-y)\,W(\d s\,\d z).
	\end{align*}
	On one hand,
	\begin{align*}
		\E(I_1^2)  &=  \int_0^t\d s\int_{x-t^{(1-\alpha)/2}}^{x+t^{(1-\alpha)/2}}\d z
			\left[G_{t-s}(z-x) - G_{t-s}(z-y)\right] ^2\\
		&\le \int_0^t\d s\int_{-\infty}^\infty\d z
			\left[G_{t-s}(z-x) - G_{t-s}(z-y)\right] ^2
			= \E\left(\left| H(t\,,x) - H(t\,,y)\right|^2 \right) \lesssim x-y;
	\end{align*}
	see \eqref{H-H}. On the other hand,
	\begin{align*}
		\E(I_2^2) &=\int_0^t\d s\int_{[x-t^{(1-\alpha)/2},\,x+t^{(1-\alpha)/2}]
			\setminus [y-t^{(1-\alpha)/2},\,y+t^{(1-\alpha)/2}]}\d z\ [G_{t-s}(z-y)]^2\\
		&=\int_0^t\d s\int_{-t^{(1-\alpha)/2}}^{x-y-t^{(1-\alpha)/2}}\d z\ [G_s(z)]^2
			+ \int_0^t\d s\int_{t^{(1-\alpha)/2}}^{x-y+t^{(1-\alpha)/2}}\d z\ [G_s(z)]^2\\
		&\le\int_0^t\frac{\d s}{\sqrt{4\pi s}}\int_{-t^{(1-\alpha)/2}}^{x-y-t^{(1-\alpha)/2}}\d z\ G_s(z)
			+ \int_0^t\frac{\d s}{\sqrt{4\pi s}}\int_{t^{(1-\alpha)/2}}^{x-y+t^{(1-\alpha)/2}}\d z\ G_s(z),
	\end{align*}
	thanks to \eqref{G}. Another appeal to \eqref{G} yields
	\[
		\int_{\pm t^{(1-\alpha)/2}}^{x-y\pm t^{(1-\alpha)/2}}G_s(z)\,\d z
		\le \frac{x-y}{\sqrt{s}}\wedge 1.
	\]
	Consequently,
	\begin{align*}
		\E(I_2^2) &\le \int_0^t\left(\frac{x-y}{\sqrt{s}}\wedge 1\right)
			\frac{\d s}{\sqrt{s}}
			=\int_0^{t\wedge |x-y|^2}\frac{\d s}{\sqrt s} +
			(x-y)\int_{t\wedge |x-y|^2}^t\frac{\d s}{s}\\
		&=2 \left\{\sqrt{t}\wedge(x-y)\right\} + (x-y)\log\left( \frac{t}{t\wedge |x-y|^2}\right)\\
		&\le 3\begin{cases}
			\sqrt t &\text{if $t\le |x-y|^2$}\\
			(x-y)\log\left(\dfrac{t}{|x-y|^2}\right)&\text{if $t>|x-y|^2$}
			\end{cases}
			\le 6(x-y)\log_+\left(\frac{1}{x-y}\right),
	\end{align*}
	uniformly for all real numbers $x>y$ and $t\in(0\,,1]$. Combine the estimates for $I_1$ and $I_2$
	to finish.
\end{proof}

The following counterpart of Lemma \ref{lem:H_alpha-H_alpha:x}
is a modulus of continuity estimate for the temporal variable of
the localized Gaussian random
field $H_\alpha$ in \eqref{H_alpha}.

\begin{lemma}\label{lem:H_alpha-H_alpha:t}
	Uniformly for every $\varepsilon\in(0\,,1)$,
	\[
		\adjustlimits\sup_{t\in[0,1]}\sup_{x\in\R}\sup_{\alpha\in(0,1)}
		\E\left(|H_\alpha(t+\varepsilon\,,x) - H_\alpha(t\,,x)|^2\right) 
		\lesssim \varepsilon^{(1-\alpha)/2}.
	\]
\end{lemma}

\begin{proof}
	Thanks to stationarity (see for example Lemma \ref{lem:H}) it suffices to consider
	only $x=0$. In that case, \eqref{H_alpha} ensures that
	we may write, for every $\alpha\in(0\,,1)$ and $t\ge0$,
	\[
		\E\left(|H_\alpha(t+\varepsilon\,,0) - H_\alpha(t\,,0)|^2\right) = I_1 + I_2 + I_3,
	\]
	where
	\begin{align*}
		I_1 &= \int_0^t\d s\int_{-t^{(1-\alpha)/2}}^{t^{(1-\alpha)/2}}\d y\,
			\left[ G_{t+\varepsilon-s}(y) - G_{t-s}(y)\right]^2,\\
		I_2 &= 2\int_0^t\d s\int_{t^{(1-\alpha)/2}}^{(t+\varepsilon)^{(1-\alpha)/2}}\d y\
			[G_{t-s+\varepsilon}(y)]^2,\\
		I_3 &= \int_t^{t+\varepsilon}\d s\int_{%
			-(t+\varepsilon)^{(1-\alpha)/2}}^{(t+\varepsilon)^{(1-\alpha)/2}}\d y\
			[G_{t-s+\varepsilon}(y)]^2.
	\end{align*}
	It is not hard to check that
	\begin{equation}\label{I1+I3}
		I_1+I_3 \le \E(|H(t+\varepsilon\,,0)-H(t\,,0)|^2) \lesssim \sqrt{\varepsilon}.
	\end{equation}
	Indeed, the inequality comes from replacing $\int_{-t^{(1-\alpha)/2}}^{t^{(1-\alpha)/2}}$
	by $\int_{-\infty}^\infty$, and the identity follows from \eqref{H-H}. 
	It therefore remains to prove that $I_2\lesssim\sqrt\varepsilon$,
	with the same parameter dependencies as in the statement of the lemma. 
	Since $[G_{s+\varepsilon}(y)]^2\le 1/(4\pi(s+\varepsilon)) 
	\lesssim \varepsilon^{-(1+\alpha)/2} s^{-(1-\alpha)/2}$
	for all $y\in\R$ and $\varepsilon,s>0$, it follows that,
	uniformly for all $t\in[0\,,1]$ and $\varepsilon>0$,
	\begin{align*}
		I_2 &\lesssim \varepsilon^{-(1+\alpha)/2}\int_0^t 
			\frac{\d s}{s^{(1-\alpha)/2} } \left[(t+\varepsilon)^{%
			(1-\alpha)/2}-t^{(1-\alpha)/2} \right]\\
			& \lesssim \varepsilon^{%
			-(1+\alpha)/2} t^{(1+\alpha)/2} \int_t^{t+\varepsilon} r^{-(1+\alpha)/2}\, \d r
		\lesssim  \varepsilon^{-(1+\alpha)/2}\,\varepsilon = \varepsilon^{(1-\alpha)/2}.
	\end{align*}
	The lemma follows from \eqref{I1+I3} and the preceding bound for $I_2$.
\end{proof}

Armed with the preceding, we can now state and prove the quantitative variation of
the main result of the earlier-announced Proposition \ref{pr:H-H_alpha:sup}.

\begin{lemma}\label{lem:H-H_alpha:prob}
	Choose and fix $\alpha \in (0\,,1)$, $R>0$, and $p \in (0\,,(1-\alpha)/4)$.
	Then, there exists $c = c(\alpha\,,R\,,p)>0$ such that
	\begin{align*}
		\P\left\{ \sup_{s \in (0,t]} \sup_{x \in [-R,R]}
		|H(s\,,x)-H_\alpha(s\,,x)| \ge \exp(-p/t^\alpha) \right\}
		\le c^{-1} \exp\left( -c\e^{c/t^\alpha}\right),
	\end{align*}
	uniformly for all $t \in (0\,,1)$.
\end{lemma}

\begin{proof}
	Throughout, we choose and fix three numbers $\alpha\in(0\,,1)$, $R>0$, and
	$p \in (0\,,(1-\alpha)/4)$.
	Recall the elementary fact that, if $X$ has a centered Gaussian distribution, then
	$\P\{|X|\ge a\} \le \exp\{- a^2/[2\E(X^2)]\}$ for all $a>0$.
	This and Lemma \ref{lem:H-H_alpha} together yield the following: There exists
	$c=c(p)>1$ such that, uniformly for all 
	$t \in (0\,,1)$ and $x\in\R$,
	\[
		\sup_{s \in (0,t]}
		\P\left\{ |H(s\,,x)-H_\alpha(s\,,x)| \ge \frac{\exp(-p/t^\alpha) }{2}\right\}
		\le \exp\left( - c \e^{c/t^\alpha} \right).
	\]
	Define
	\begin{equation}\label{X_L}
		\mathbb{X}_L = \bigcup_{j=0}^{L^2}\bigcup_{k=-L}^L
		\left\{\left( \frac{j}{L^2}\,,\frac{k}{L}\right)\right\}\qquad\forall L\in\N,
	\end{equation}
	and deduce from the preceding that
	\begin{align*}
		&\P\left\{ \max_{(s,x)\in\mathbb{X}_L\cap((0,t]\times[-R,R])}
			|H(s\,,x)-H_\alpha(s\,,x)| \ge \frac{\exp(-p/t^\alpha) }{2}\right\}\\
		&\hskip1in\le (L^2+1)(1+2L) \exp\left( -c\e^{c/t^\alpha} \right)
			\le 6 L^3\exp\left( -c\e^{c/t^\alpha} \right).
	\end{align*}
	Therefore,
	\[
		\P\left\{ \sup_{(s,x)\in (0,t] \times[-R,R]}
		|H(s\,,x)-H_\alpha(s\,,x)| \ge \exp(-p/t^\alpha)\right\}
		\le 6 L^3\exp\left( -c\e^{c/t^\alpha} \right)
		+ \mathscr{P}_1 + \mathscr{P}_2,
	\]
	where
	\begin{align*}
		\mathscr{P}_1 &= \P\left\{ \sup_{(s,x)\in (0,t] \times[-R,R]}
    			\sup_{\substack{r:|r-s|\le 1/L^2\\
			y:|x-y|\le 1/L}} |H(s\,,x)-H(r\,,y)| \ge \frac{\exp(-p/t^\alpha)}{4}\right\},\\
		\mathscr{P}_2 &=P\left\{ \sup_{(s,x)\in (0,t] \times[-R,R]}
    			\sup_{\substack{r:|r-s|\le 1/L^2\\
			y:|x-y|\le 1/L}} |H_\alpha(s\,,x)-H_\alpha(r\,,y)| \ge \frac{\exp(-p/t^\alpha)}{4}\right\},
	\end{align*}
	and where the dependencies on the parameters $(n\,,R\,,\varepsilon\,,\ldots)$ are supressed 
	to simplify the exposition and the notation.
	For $\gamma \in (0\,, (1-\alpha)/4)$, 
	Lemmas \ref{lem:H_alpha-H_alpha:x} and \ref{lem:H_alpha-H_alpha:t} can be combined with
	a standard metric entropy argument (see, e.g., \cite[Proposition A.1]{DKN07})
	to yield the following:
	\[
		\E\left(\sup_{\substack{0\le s<t\le1\\
		-R\le x<y\le R}} \frac{|H_\alpha(t\,,x)-H_\alpha(s\,,y)|}{|t-s|^\gamma+|x-y|^{2\gamma}}\right)<\infty.
	\]
	And, as is well known, the same inequality holds (for more or less the same type of reasons)
	when we replace $H_\alpha$ by $H$ everywhere. 
	Therefore, the concentration properties of Gaussian measure (see  Ledoux \cite{Ledoux}) yields
	a constant $\zeta=\zeta(\alpha\,,R\,,p\,,\gamma)>0$ such that
	\begin{equation}\label{A1A2}\begin{split}
		A_1&=\E\left[\sup_{\substack{0\le s<t\le1\\
			-R\le x<y\le R}} \exp\left(\zeta\frac{|H(t\,,x)-H(s\,,y)|^2}{
			|t-s|^{2p}+|x-y|^{4p}}\right)\right]<\infty,
			\text{ and}\\
		A_2&=\E\left[\sup_{\substack{0\le s<t\le1\\
			-R\le x<y\le R}} \exp\left(\zeta
			\frac{|H_\alpha(t\,,x)-H_\alpha(s\,,y)|^2}{
			|t-s|^{2p}+|x-y|^{4p}}\right)\right]<\infty.
	\end{split}\end{equation}
	The preceding and Chebyshev's inequality together yield
	\[
		\mathscr{P}_i \le A_i\exp\left( - \frac{\zeta L^{4\gamma}
		\exp(-2p/t^\alpha) }{16}\right)
		\qquad\forall i=1,2.
	\]
	Consequently,
	\begin{align*}
		&\P\left\{ \sup_{(s,x)\in (0,t] \times[-R,R]}
		|H(s\,,x)-H_\alpha(s\,,x)|\ge \e^{-p/t^\alpha} \right\}\\
		&\hskip1in\le 6 L^3\exp\left( - c\, \e^{c/t^\alpha}\right) + (A_1+A_2)
		\exp\left( - c L^{4\gamma} \e^{-2p/t^\alpha}\right).
	\end{align*}
	Finally, we may take $L = \e^{A/t^\alpha}$, where $A>p/(2\gamma)$, in order to deduce that
	\begin{align*}
		&\P\left\{ \sup_{(s,x)\in (0,t] \times[-R,R]}
			|H(s\,,x)-H_\alpha(s\,,x)|\ge \exp(-p/t^\alpha) \right\}\\
		&\le 6 \exp\left( \frac{A}{t^\alpha} - c\, \e^{c/t^\alpha}\right) + (A_1+A_2)
			\exp\left( - c \e^{(4\gamma A-2p)/t^\alpha}\right)
			\le c^{-1} \exp\left( -c \e^{c/t^\alpha}\right).
	\end{align*}
	This completes the proof of Lemma \ref{lem:H-H_alpha:prob}.
\end{proof}

\begin{proof}[Proof of Proposition \ref{pr:H-H_alpha:sup}]
	It follows readily from Lemma \ref{lem:H-H_alpha:prob} that
	\[
		\sum_{n=1}^\infty \P\left\{ \sup_{(t,x)\in [\frac{1}{n+1},\frac{1}{n}]\times[-R,R]}
		\frac{|H(t\,,x)-H_\alpha(t\,,x)|}{t^{1/4}} \ge  \varepsilon \right\}<\infty.
	\]
	Proposition \ref{pr:H-H_alpha:sup} follows from the above and the Borel-Cantelli lemma.
\end{proof}

\begin{lemma}\label{lem:H_alpha}
	Choose and fix $\alpha \in (0\,,1)$, $\theta>0$ and $q>1$. Then,
	\[
		\P\left\{ |H_\alpha(s\,,0)|\le \theta s^{1/4}\,\,\forall s\in [a\,,b]\right\}
		=  (a/b)^{\lambda(\theta) + \mathscr{o}(1)}
	\]
	as $\max\{b\,,a/b\} \to 0^+$, subject to $b^q \le a \le b$.
\end{lemma}

\begin{proof}
	By disjointness of integration domains and the properties of Wiener 
	integrals, $H - H_\alpha$ is a centered Gaussian random field that is independent
	of $H_\alpha$ (see \eqref{mild:H} and \eqref{H_alpha}).
	Anderson's inequality \cite{Anderson} implies that for any $0<a<b$ and $f \in C([a\,,b])$,
	\[
		\P\left\{ |H_\alpha(s\,,0) + f(s)| \le \theta s^{1/4} \, 
		\forall s \in [a\,,b] \right\} \le \P\left\{ |H_\alpha(s\,,0)| 
		\le \theta s^{1/4} \, \forall s \in [a\,,b] \right\}.
	\]
	By conditioning on $H-H_\alpha$ and applying the above inequality, we deduce that
	\[
		\P\left\{ |H(s\,,0)|\le \theta s^{1/4}\,\,\forall s\in[a\,,b]\right\}
		\le\P\left\{ |H_\alpha(s\,,0)|\le \theta s^{1/4}\,\,\forall s\in[a\,,b]\right\}.
	\]
	Thanks to the scaling properties of $H$, the left-hand side is equal to
	\[
		\P\left\{ |H(s\,,0)|\le \theta s^{1/4}\,\,\forall s\in[a/b\,,1]\right\}
		=(a/b)^{\lambda(\theta)+\mathscr{o}(1)},
	\]
	uniformly for all $0<a<b$ such that $a/b\to0^+$. This proves that
	\begin{equation}\label{HHLB}
		\P\left\{ |H_\alpha(s\,,0)|\le \theta s^{1/4}\,\,\forall s\in[a/b\,,1]\right\}
		\ge (a/b)^{\lambda(\theta)+\mathscr{o}(1)},
	\end{equation}
	uniformly for all $0<a<b$ such that $a/b\to0^+$; this is a slightly stronger
	lower bound than what we need since $a$ and $b$ do not need to converge to zero 
	individually for this bound, only the ratio needs to. For the complementary bound
	we may observe that for every $\delta\in(0\,,1)$
	there exists $b_0 = b_0(\alpha\,,\delta\,,\theta)\in(0\,,1)$ such that
	\begin{align*}
		&\P\left\{ |H_\alpha(s\,,0)|\le \theta s^{1/4}\,\,\forall s\in[a\,,b]\right\}\\
			&\le \P\left\{ |H(s\,,0)|\le \theta s^{1/4} + 
			\e^{-\beta/s^\alpha}\,\,\forall s\in[a\,,b]\right\}
			+ c^{-1}\exp\left(-c\e^{c/b^\alpha}\right)\\
		&\le \P\left\{ |H(s\,,0)|\le (\theta+\delta) s^{1/4}\,\,\forall s\in[a\,,b]\right\}
			+ c^{-1}\exp\left(-c\e^{c/b^\alpha}\right),
	\end{align*}
	uniformly for all $0<a<b<b_0$. In particular, the scaling properties of
	$H$ imply that, uniformly for all $0<a<b<b_0$,
	\begin{align*}
		&\P\left\{ |H_\alpha(s\,,0)|\le \theta s^{1/4}\,\,\forall s\in[a\,,b]\right\}\\
		&\hskip1in
			\le (a/b)^{\lambda(\theta+\delta)+\mathscr{o}(1)}
			+ c^{-1}\exp\left(-c\e^{c/b^\alpha}\right)= (a/b)^{\lambda(\theta+\delta)+\mathscr{o}(1)},
	\end{align*}
	as $\max\{b\,,a/b\}\to0^+$, subject to $b^q \le a \le b$.
	In other words, together with \eqref{HHLB}, this proves that
	whenever $b^q \le a \le b$ and $\max\{b\,,a/b\}\to0^+$,
	\begin{align*}
		\lambda(\theta+\delta)+\mathscr{o}(1) \le 
		\frac{\log\P\left\{ |H_\alpha(s\,,0)|\le \theta s^{1/4}\,\,\forall s\in[a\,,b]\right\}}{\log (b/a)}
		\le \lambda(\theta)+\mathscr{o}(1),
	\end{align*}
	for every fixed choice of
	$\delta\in(0\,,1)$, as well as $\theta>0$. 
	Since $\lambda$ is continuous (see Proposition \ref{pr:BC:H}),
	we may let $\delta$ tend to zero in order to obtain the result.
\end{proof}

\section{Proof of Theorem \ref{th:BC:u}, Part (2)}

We begin with the presentation of Part (2) of Theorem \ref{th:BC:u}.
Before we do that, let us observe that, thanks to Proposition \ref{pr:L_k:supE}
and the Borel-Cantelli lemma,
\[
	\P\left\{ \lim_{t\to0^+} \sup_{x \in K}\frac{|\mathscr{E}(t\,,x)|}{t^{1/4}} =0
	\quad\forall x\in\R\right\}=1,
\]
where the space-time random field $\mathscr{E}$ was defined in \eqref{E(t,x)}.
In particular, it follows that the slow point of $u$, as defined by \eqref{S_u(theta)}
can  be written in terms of the Gaussian random field $H$ of \eqref{H} as 
follows: With probability one,
\begin{equation}\label{S_H}
	\fS(\theta) = \left\{ x\in\R:\ 
	\limsup_{t\to0^+}\frac{|H(t\,,x)|}{t^{1/4}} \le \theta \right\}
	\quad\forall\theta>0.
\end{equation}
To be precise, let us write the right-hand side as
$\fS_H(\theta)$. Then, we have shown that the
set $\cup_{\theta>0}\{\fS(\theta)\neq\fS_H(\theta)\}$ is a subset of
a $\P$-null set, and hence is $\P$-null thanks to the fact that our probability space is
complete.
In this way, we are permitted to use \eqref{S_H} as a definition of $\fS(\theta)$,
which we will. This reduces our problem to one about the Gaussian random field $H$.
With the preceding under way, the real work can now begin.

\begin{theorem}\label{th:H}
	Choose and fix a nonrandom compact set
	$K \subset \R$.
	Then there is a $\P$-null set off which
	\begin{align*}
		\lambda^{-1}\left(\tfrac12 \dimm K \right) \le 
		\adjustlimits \inf_{x \in K} \limsup_{t \to 0^+} \frac{|H(t\,,x)|}{t^{1/4}}
		\le \lambda^{-1}\left( \tfrac12 \dimh K \right),
	\end{align*}
	where $\lambda^{-1}(0) := \inf_{\theta>0}\lambda^{-1}(\theta)=\infty$.
\end{theorem}

This implies Theorem \ref{th:BC:u}(2).
The proof of Theorem \ref{th:H} itself is divided in two parts and will be 
given in the \S\S\ref{s:th:H:lb} and \ref{s:th:H:ub} that follow.

\subsection{Proof of Theorem \ref{th:H}: Lower bound}\label{s:th:H:lb}
Let $K\subset\R$ be a nonempty compact set and choose and fix
a number $\varepsilon>0$. Recall that a finite set $L\subset K$ is said to be an
\emph{$\varepsilon$-covering}
of $K$ if every $x\in K$ is within $\varepsilon$ of some $y\in L$. 
Throughout, the number 
$\mathscr{N}_\varepsilon(K)$ will denote the cardinality of the smallest $\varepsilon$-covering of $K$.
Because $K$ is totally bounded, $\mathscr{N}_\varepsilon(K)<\infty$ for every $\varepsilon>0$.
The \emph{upper Minkowski dimension of $K$} is the number
\begin{equation}\label{dimm1}
	\dimm K = \limsup_{\varepsilon\to0+}
	\frac{\log \mathscr{N}_\varepsilon(K)}{|\log\varepsilon|}.
\end{equation}
Some authors refer to $\dimm$ alternatively as the  \emph{upper box dimension};
see Falconer \cite{Falconer}.

Since $K$ is compact, it is contained in a bounded interval $[a\,,b]$ that we fix.
Choose and fix numbers $\theta,\eta,\beta > 0$ and $\rho\in (0\,,1)$ such that
\begin{align}\label{lambda>}
	\lambda(\theta+\eta) > \frac{(1+\rho)\beta}{2}
	> \frac{1+\rho}{2} \dimm K.
\end{align}
This can always be done since $\lim_{c\to0^+}\lambda(c)=\infty$
[Proposition \ref{pr:BC:H}].
Since $\beta> \dimm K$, it follows from the definition of upper Minkowski dimension
that there exists $\varepsilon_0 \in (0\,,1)$ such that
\begin{align}\label{packing:UB}
	\mathscr{N}_{\varepsilon}(K) \le \varepsilon^{-\beta} \qquad \forall \varepsilon \in (0\,,\varepsilon_0].
\end{align}
Choose and fix $\alpha \in (0\,,1)$ such that
\begin{align}\label{rho:alpha:>1}
	(1+\rho)(1-\alpha)>1.
\end{align}
In light of \eqref{lambda>}, we choose and fix $\mu>1$ such that
\begin{align}\label{mu}
	\mu\left(\lambda(\theta+\eta)-\frac{(1+\rho)\beta}{2}\right)
	> \lambda(\theta+\eta) + \frac{1-\alpha}{2},
\end{align}
and define
\begin{align}\label{eps_n}
	\varepsilon_n = \exp(-\mu^n) \qquad \forall n\in\N.
\end{align}
Thanks to \eqref{packing:UB}, for every large-enough $n \in \N$ there exists a
finite set $F_n\subset K$ such that the following properties hold:
\begin{enumerate}[\bf Property A.]
	\item For every $x \in K$, there exists $y \in F_n$ such that 
		$|x-y| \le \varepsilon_{n+1}^{(1+\rho)/2}$;	\label{pack}
	\item Whenever $w\,,z \in F_n$ are distinct, 
		$|w-z| \ge 2\varepsilon_{n+1}^{(1+\rho)/2}$; and
	\item	$|F_n| \lesssim \varepsilon_{n+1}^{-(1+\rho)\beta/2}$ for all $n\in\N$,
		where $|\cdots|$ denotes cardinality.
	\label{|F|}
\end{enumerate}
For every $n \gg 1$, we consider a uniform partition of the interval $[a\,,b]$ 
as follows:
\[
	[a\,,b] = [a_{n,1}\,,a_{n,2}] \cup \cdots \cup [a_{n,N}\,,a_{n,N+1}],	
\]
where
\begin{align}\label{N}
	N=N(n) = \left\lfloor \frac{b-a}{2\varepsilon_n^{(1-\alpha)/2}} 
	\right\rfloor \asymp \varepsilon_n^{-(1-\alpha)/2}
	\quad \text{and}\quad
	a_{n,i+1}-a_{n,i} = \frac{b-a}{N} \ge 2 \varepsilon_n^{(1-\alpha)/2}.
\end{align}
For every $1 \le i \le N$, let
\begin{align*}
	K_i = K_{n,i} = K \cap [a_{n,i}\,,a_{n,i+1}].
\end{align*}
Owing to Property \ref{|F|}, we have
\begin{align}\label{FnK}
	|F_n \cap K_i| \le |F_n|
	\lesssim \varepsilon_{n+1}^{-(1+\rho)\beta/2} \qquad (1 \le i \le N).
\end{align}
We may write $K = I \cup J$, where
\[
	I =  K_{1} \cup K_{3} \cup \ldots \quad \text{and} \quad
	J = K_{2} \cup K_{4} \cup \ldots\,.
\]
Recall the random field $H_\alpha$ from \eqref{H_alpha}.
By Boole's inequality,
\begin{align*}
	&\P\left\{ \adjustlimits\inf_{x \in K} 
		\sup_{t \in [\varepsilon_{n+1},\varepsilon_n]}
		\frac{|H_\alpha(t\,,x)|}{t^{1/4}} \le \theta \right\}\\
	&\le \P\left\{ \adjustlimits
		\inf_{x \in I} \sup_{t \in [\varepsilon_{n+1},\varepsilon_n]}
		\frac{|H_\alpha(t\,,x)|}{t^{1/4}} \le \theta \right\} + 
		\P\left\{ \adjustlimits
		\inf_{x \in J} \sup_{t \in [\varepsilon_{n+1},\varepsilon_n]}
		\frac{|H_\alpha(t\,,x)|}{t^{1/4}} \le \theta \right\}.
\end{align*}
Thanks to the second assertion in \eqref{N}, 
the space-time random fields
\[
	\{H_\alpha(t\,,x)\}_{t \in [\varepsilon_{n+1},\varepsilon_n], x \in K_i}
	\quad\text{and}\quad
	\{H_\alpha(t\,,x)\}_{t \in [\varepsilon_{n+1},\varepsilon_n], x \in K_{i+2}}
\]
are independent from one another. Consequently,
\begin{align*}
	&\P\left\{ \adjustlimits
		\inf_{x \in I} \sup_{t \in [\varepsilon_{n+1},\varepsilon_n]} 
		\frac{|H_\alpha(t\,,x)|}{t^{1/4}} \le \theta \right\}
		= 1-\P\left\{ \adjustlimits
		\inf_{x \in I} \sup_{t \in [\varepsilon_{n+1},\varepsilon_n]}
		\frac{|H_\alpha(t\,,x)|}{t^{1/4}} > \theta \right\}\\
	&= 1 - \prod_{i: K_i \subset I} \P\left\{ \adjustlimits\inf_{x \in K_i}
		\sup_{t \in [\varepsilon_{n+1},\varepsilon_n]}
			\frac{|H_\alpha(t\,,x)|}{t^{1/4}} > \theta \right\}\\
	&= 1 - \prod_{i\in\N: K_i \subset I}\left[ 1 - \P\left\{ \adjustlimits
		\inf_{x \in K_i} \sup_{t \in [\varepsilon_{n+1},\varepsilon_n]}
		\frac{|H_\alpha(t\,,x)|}{t^{1/4}} \le \theta \right\} \right].
\end{align*}
Consider an index $i\in\N$ such that $K_i \subset I$. 
Thanks to Property \ref{pack}, we may use interpolation 
to deduce that
\begin{align*}
	&\P\left\{ \adjustlimits\inf_{x \in K_i}
		\sup_{t \in [\varepsilon_{n+1},\varepsilon_n]}
		\frac{|H_\alpha(t\,,x)|}{t^{1/4}} \le \theta \right\}
		\le \P\left\{ \adjustlimits
		\min_{y \in F_n \cap K_i} \sup_{t \in [\varepsilon_{n+1},\varepsilon_n]}
		\frac{|H_\alpha(t\,,y)|}{t^{1/4}} \le \theta + \eta \right\}\\
	&\hskip1.5in + \P\left\{ \sup_{\substack{x,y \in [a,b]:\\ |x-y| \le
		\varepsilon_{n+1}^{(1+\rho)/2}}}
		\sup_{t \in [\varepsilon_{n+1},\varepsilon_n]} 
		\frac{|H_\alpha(t\,,x) - H_\alpha(t\,,y)|}{t^{1/4}} > \eta \right\}.
\end{align*}
A union bound (Boole's inequality), \eqref{FnK}, and Lemma \ref{lem:H_alpha}, 
together yield
\begin{align*}
	&\P\left\{ \adjustlimits\min_{y \in F_n \cap K_i} \sup_{t \in [\varepsilon_{n+1},\varepsilon_n]}
		\frac{|H_\alpha(t\,,y)|}{t^{1/4}} \le \theta + \eta \right\}\\
	& \le |F_n \cap K_i| \sup_{y \in K} \P\left\{ 
		\sup_{t \in [\varepsilon_{n+1},\varepsilon_n]} 
		\frac{|H_\alpha(t\,,y)|}{t^{1/4}} \le \theta + \eta \right\}\\
	&\lesssim \varepsilon_{n+1}^{-\frac{(1+\rho)\beta}{2}}
		(\varepsilon_{n+1}/\varepsilon_n)^{\lambda(\theta+\eta) + \mathscr{o}(1)}
		= \varepsilon_{n+1}^{\lambda(\theta+\eta)-\frac{(1+\rho)\beta}{2} 
		+ \mathscr{o}(1)} \varepsilon_n^{-\lambda(\theta+\eta) + \mathscr{o}(1)},
\end{align*}
as $n \to \infty$.
Choose and fix $p \in (0\,,(1-\alpha)/4)$ such that $\delta :=  (1+\rho)p -1/4>0$, 
possible provided that $p$ is sufficiently close to $(1-\alpha)/4$ due to \eqref{rho:alpha:>1}.
Thanks to \eqref{A1A2}, there exists a constant $c>0$ such that
\begin{align*}
	&\P\left\{ \sup_{\substack{x,y \in [a,b]:\\ |x-y| \le \varepsilon_{n+1}^{(1+\rho)/2}}}
		\sup_{t \in [\varepsilon_{n+1},\varepsilon_n]}
		\frac{|H_\alpha(t\,,x) - H_\alpha(t\,,y)|}{t^{1/4}} > \eta \right\}\\
	&\le \P\left\{ \sup_{x,y \in [a,b]} \sup_{t \in [0,1]}
		\frac{|H_\alpha(t\,,x) - H_\alpha(t\,,y)|}{|x-y|^{2p}} > \eta 
		\varepsilon_{n+1}^{\frac14 - (1+\rho)p}\right\}
		\lesssim \exp\left( - \frac{\eta^2}{c \varepsilon_{n+1}^{2\delta}} \right),
\end{align*}
uniformly for all $n \in \N$.
It follows that
\begin{align*}
	&\P\left\{ \adjustlimits
		\inf_{x \in K} \sup_{t \in [\varepsilon_{n+1},\varepsilon_n]}
		\frac{|H_\alpha(t\,,x)|}{t^{1/4}} \le \theta \right\}\\
	&\lesssim \varepsilon_{n+1}^{\lambda(\theta+\eta)-\frac{(1+\rho)\beta}{2}
		+ \mathscr{o}(1)} \varepsilon_n^{-\lambda(\theta+\eta) + \mathscr{o}(1)}
		 + \exp\left(- \frac{\eta^2 }{c\varepsilon_{n+1}^{2\delta}}\right)
		\quad\text{as $n \to \infty$.}
\end{align*}
The same estimate applies to every index $i\in\N$ that satisfies either
$K_i \subset I$ or $K_i \subset J$.
Therefore, we may combine the above estimates above, and apply \eqref{N}, 
in order to see that there exists a constant $C>0$ such that, as $n\to\infty$,
\begin{align*}
	&\P\left\{\adjustlimits
		\inf_{x \in K} \sup_{t \in [\varepsilon_{n+1},\varepsilon_n]}
		\frac{|H_\alpha(t\,,x)|}{t^{1/4}} \le \theta \right\}\\
	& \le 2 \left[ 1 - \left( 1 - C 
		\varepsilon_{n+1}^{\lambda(\theta+\eta)-\frac{(1+\rho)\beta}{2} 
		+ \mathscr{o}(1)} \varepsilon_n^{-\lambda(\theta+\eta) + \mathscr{o}(1)}
	- C\exp\left\{- \frac{\eta^2 }{C\varepsilon_{n+1}^{2\delta}}\right\}
		\right)^{C \varepsilon_n^{-\frac{1-\alpha}{2}}} \right].
\end{align*}
According to \eqref{mu} and \eqref{eps_n}, as $n\to\infty$,
\begin{align*}
	&\varepsilon_{n+1}^{\lambda(\theta+\eta)-\frac{(1+\rho)\beta}{2} 
		+ \mathscr{o}(1)} \varepsilon_n^{-\lambda(\theta+\eta) + \mathscr{o}(1)}
		= \varepsilon_n^{\mu(\lambda(\theta+\eta)-\frac{(1+\rho)\beta}{2})
		-\lambda(\theta+\eta) + \mathscr{o}(1)} = \mathscr{o}(1).
\end{align*}
Therefore, the elementary inequality $\exp(-z) \le 1-(z/2)$,
valid for all $z \in [0\,,1]$, \eqref{mu}, and \eqref{eps_n} together yield
\begin{align*}
	&\P\left\{ \adjustlimits
		\inf_{x \in K} \sup_{t \in [\varepsilon_{n+1},\varepsilon_n]}
		\frac{|H_\alpha(t\,,x)|}{t^{1/4}} \le \theta \right\}\\
	&\le2\left[ 1 - \exp\left\{ -2C^2 \varepsilon_n^{-\frac{1-\alpha}{2}}
		\left(\varepsilon_{n}^{\mu(\lambda(\theta+\eta)-
		\frac{(1+\rho)\beta}{2})-\lambda(\theta+\eta)
		+\mathscr{o}(1)} + \exp\left[-\frac{\eta^2}{C\varepsilon_{n+1}^{2\delta}}\right]
		\right)\right\} \right]\\
	&= 2\left[ 1 - \exp\left\{ -2C^2 \left(\e^{ - \mu^n 
		[ \mu ( \lambda(\theta+\eta) - \frac{(1+\rho)\beta}{2} ) -\frac{1-\alpha}{2} - 
		\lambda(\theta+\eta) + \mathscr{o}(1) ]} + 
		\varepsilon_n^{-\frac{1-\alpha}{2}} 
		\exp\left[-\frac{\eta^2}{C\varepsilon_{n+1}^{2\delta}}\right]
		\right)\right\} \right]\\
	&=\mathscr{o}(1) \quad \text{as $n\to\infty$.}
\end{align*}
Let $A_n$ denote the complement of the event in the first line of the last display.
It follows that $\liminf_{n\to\infty}\P(A_n^c)=0.$
Because $\P\{A_n\text{\ \rm i.o.}\} \ge 1-\liminf_{n\to\infty}\P(A_n^c),$
it follows that, almost surely, infinitely many of the events $A_n$ must occur. 
Thus, we can deduce from Proposition \ref{pr:H-H_alpha:sup} that
\begin{align}\label{E:inflim:LB}
	\adjustlimits\inf_{x\in K}\limsup_{t \to 0^+} \frac{|H(t\,,x)|}{t^{1/4}} 
	= \adjustlimits\inf_{x\in K}\limsup_{t \to 0^+} 
	\frac{|H_\alpha(t\,,x)|}{t^{1/4}} \ge \theta\quad\text{a.s.}
\end{align}
Recall that $\theta$, $\eta$, $\rho$, $\beta$ satisfy \eqref{lambda>}.
Because $\lambda$ is strictly decreasing and continuous 
(see Proposition \ref{pr:BC:H}), we can let $\eta \downarrow 0$, 
$\rho \downarrow 0$, and then $\beta \downarrow \dimm K$, all along rationals, in order
to see that
\eqref{E:inflim:LB} holds for all $\theta>0$ such that $\lambda(\theta)>\tfrac12\dimm K$.
Take supremum over all such $\theta>0$ in order to find that
\[
	\adjustlimits\inf_{x\in K}\limsup_{t \to 0^+} 
	\frac{|H(t\,,x)|}{t^{1/4}} \ge \sup\{\theta>0: 
	\lambda(\theta) > \tfrac12 \dimm K\}
	\quad\text{a.s.}
\]
Thanks to the monotonicity of $\lambda$ (Proposition \ref{pr:BC:H}),
\[
	\sup\{\theta>0: \lambda(\theta) > \tfrac12 \dimm K\}
	=\inf\{\theta>0: \lambda(\theta) < \tfrac12 \dimm K\},
\]
and both of these quantities are equal to $\lambda^{-1}(\frac12\dimm K)$,
where $\lambda^{-1}(0)=\infty$.
This completes the proof of the lower bound for Theorem \ref{th:H}.
\qed

\begin{remark}
	Property \ref{|F|}, crucial in the proof,
	relied on the use of the upper Minkowski dimension
	(and not the Hausdorff dimension).
	This is why our method does not yield a lower bound with $\dimh K$.
\end{remark}

\subsection{Proof of Theorem \ref{th:H}: Upper bound}\label{s:th:H:ub}

We now turn to the proof of the upper bound.

\begin{proof}
	We may, and will assume without loss of generality that $\dimh K>0$
	for there is nothing to prove otherwise. Let us fix two numbers
	$\beta,\theta>0$ such that
	\[
		2\lambda(\theta) < \beta < \dimh K.
	\]
	By Frostman's lemma (see Falconer \cite{Falconer}), there exists 
	a probability measure $\mu$ on $K$ such that
	\begin{align}\label{Frostman:meas}
		\adjustlimits \sup_{r>0}\sup_{x\in\R}
		\frac{\mu([x-r\,, x+r])}{r^\beta}<\infty.
	\end{align}
	Fix $\alpha \in (0\,,1)$  such that
	\begin{align}\label{alpha:rho}
		\frac{\lambda(\theta)}{\beta} < \frac{1-\alpha}{2} < \frac12,
		\qquad\text{set}\quad \rho = \frac{1-\alpha}{2},
	\end{align}
	and recall the localization $\{H_\alpha(t\,,x)\}_{t>0,x \in \R}$ 
	of $H$ from \eqref{H_alpha}.
	For every $\varepsilon_2>\varepsilon_1 > 0$ and $x \in \R$, define the event
	\[
		\mathcal{E}(\varepsilon_1\,,\varepsilon_2\,,x)  =  
		\left\{  \omega\in\Omega:\
		\sup_{t \in [\varepsilon_1,\varepsilon_2]} 
		\frac{|H_\alpha(t\,,x)|(\omega)}{t^{1/4}} \le \theta \right\}.
	\]
	Choose and fix some $q>1$.
	The spatial stationarity of $H_\alpha$ and Lemma \ref{lem:H_alpha} together
	imply that
	\begin{align}\label{varphi}
		\P(\mathcal{E}(\varepsilon_1\,,\varepsilon_2\,,x))
		= \P(\mathcal{E}(\varepsilon_1\,,\varepsilon_2\,,0)) = 
		(\varepsilon_1/\varepsilon_2)^{\lambda(\theta)+\mathscr{o}(1)},
	\end{align}
	as $\max\{\varepsilon_2\,,\varepsilon_1/\varepsilon_2\} \to 0^+$, 
	subject to $\varepsilon_2^q \le \varepsilon_1 \le \varepsilon_2$, uniformly for all $x \in K$.
	For every $\varepsilon_2>\varepsilon_1>0$, consider the random variable
	\[
		X_{\varepsilon_1,\varepsilon_2}  =  
		\int_K \frac{\1_{\mathcal{E}(\varepsilon_1,\varepsilon_2,x)}}{%
		\P(\mathcal{E}(\varepsilon_1,\varepsilon_2,x))} \mu( \d x).
	\]
	Define $h_m = \exp(-q^m)$ for every $m \in \N$. 
	We aim to prove that 
	\[
		\sum_{m=1}^\infty \P\{ X_{h_{m+1}, h_m} = 0\} < \infty
	\]
	for a suitable choice of $q>1$.
	By the first identity in \eqref{varphi},
	\[
		\E(X_{\varepsilon_1,\varepsilon_2}^2)
		= \iint_{K\times K} \frac{\P(\mathcal{E}(%
		\varepsilon_1\,,\varepsilon_2\,,x) \cap \mathcal{E}(%
		\varepsilon_1\,,\varepsilon_2\,,y))}{\P(\mathcal{E}(%
		\varepsilon_1\,,\varepsilon_2\,,0))^2}\, \mu(\d x)\, \mu(\d y).
	\]
	We may decompose $\E(X_{\varepsilon_1,\varepsilon_2}^2)$ as follows:
	\[
		\E(X_{\varepsilon_1,\varepsilon_2}^2)
		= I_1(\varepsilon_1\,,\varepsilon_2)
		+I_2(\varepsilon_1\,,\varepsilon_2)
		+I_3(\varepsilon_1\,,\varepsilon_2),
	\]
	where
	\begin{align*}
		&I_1 = I_1(\varepsilon_1\,,\varepsilon_2) = 
			\iint_{\substack{x,y \in K:\\|x-y|\le 2\varepsilon_1^\rho}} 
			\frac{\P(\mathcal{E}(\varepsilon_1\,,\varepsilon_2\,,x) \cap 
			\mathcal{E}(\varepsilon_1\,,\varepsilon_2\,,y))}{%
			\P(\mathcal{E}(\varepsilon_1\,,\varepsilon_2\,,0))^2}\, 
			\mu(\d x)\, \mu(\d y),\\
		&I_2 = I_2(\varepsilon_1\,,\varepsilon_2) = 
			\iint_{\substack{x,y \in K:\\ 2\varepsilon_1^\rho < |x-y|\le 
			2\varepsilon_2^\rho}} \frac{\P(\mathcal{E}(\varepsilon_1\,,\varepsilon_2\,,x)
			\cap \mathcal{E}(\varepsilon_1\,,\varepsilon_2\,,y))}{%
			\P(\mathcal{E}(\varepsilon_1\,,\varepsilon_2\,,0))^2}\, \mu(\d x)\, \mu(\d y),\\
		&I_3 = I_3(\varepsilon_1\,,\varepsilon_2) = \iint_{%
			\substack{x,y \in K:\\|x-y|>2\varepsilon_2^\rho}} 
			\frac{\P(\mathcal{E}(\varepsilon_1\,,\varepsilon_2\,,x) 
			\cap \mathcal{E}(\varepsilon_1\,,\varepsilon_2\,,y))}{%
			\P(\mathcal{E}(\varepsilon_1\,,\varepsilon_2\,,0))^2}\, \mu(\d x)\, \mu(\d y).
	\end{align*}
	The trivial bound $\P(\mathcal{E}_1\cap\mathcal{E}_2) \le \P(\mathcal{E}_1)$ 
	is valid for all events $\mathcal{E}_1$ and $\mathcal{E}_2$. 
	This, together with \eqref{Frostman:meas}, \eqref{alpha:rho} and \eqref{varphi},  yields
	\[
		I_1 \lesssim \frac{\varepsilon_1^{\beta\rho}}{\P(\mathcal{E}(
		\varepsilon_1\,,\varepsilon_2\,,0))} = \varepsilon_1^{%
		\beta\rho-\lambda(\theta)+\mathscr{o}(1)} \varepsilon_2^{\lambda(\theta)
		+\mathscr{o}(1)} = \mathscr{o}(1)
	\]
	as $\max\{\varepsilon_2\,,\varepsilon_1/\varepsilon_2\} 
	\to 0^+$, subject to $\varepsilon_2^\rho \le \varepsilon_1 \le \varepsilon_2$.
	
	Next, if $|x-y|>2\varepsilon_2^\rho$ then,
	by the disjointness of integration domains, the processes
	\[
		\{ H_\alpha(t\,,x) \}_{t \in [\varepsilon_1,\varepsilon_2]} \quad \text{and} \quad
		\{ H_\alpha(t\,,y) \}_{t \in [\varepsilon_1,\varepsilon_2]}
	\]
	are independent (see \eqref{H_alpha} and \eqref{alpha:rho}).
	Hence, the events
	$\mathcal{E}(\varepsilon_1\,,\varepsilon_2\,,x)$
	and $\mathcal{E}(\varepsilon_1\,,\varepsilon_2\,,y)$ are independent. Therefore, $I_3 \le 1$.
	
	We estimate $I_2$ as follows: 
	If $2\varepsilon_1^\rho < |x-y| \le 2\varepsilon_2^\rho$, then 
	\begin{align*}
		&\P\left\{ \sup_{t \in [\varepsilon_1,\varepsilon_2]} 
			\frac{|H_\alpha(t\,,x)|}{t^{1/4}} \le \theta\ , \, 
			\sup_{t \in [\varepsilon_1,\varepsilon_2]} 
			\frac{|H_\alpha(t\,,y)|}{t^{1/4}} \le \theta \right\}\\
		& \le \P\left\{ \sup_{t \in [\varepsilon_1, (|x-y|/2)^{1/\rho}]} 
			\frac{|H_\alpha(t\,,x)|}{t^{1/4}} \le \theta \ , \, 
			\sup_{t\in [\varepsilon_1,(|x-y|/2)^{1/\rho}]} 
			\frac{|H_\alpha(t\,,y)|}{t^{1/4}} \le \theta \right\}\\
		&= \P\left\{ \sup_{t \in [\varepsilon_1, (|x-y|/2)^{1/\rho}]} 
			\frac{|H_\alpha(t\,,x)|}{t^{1/4}} \le \theta \right\} 
			\P\left\{ \sup_{t\in [\varepsilon_1,(|x-y|/2)^{1/\rho}]} 
			\frac{|H_\alpha(t\,,y)|}{t^{1/4}} \le \theta \right\},
	\end{align*}
	where the last line follows from independence of the two events, 
	due to independence of the processes
	\[
		\{ H_\alpha(t\,,x) \}_{t \in [\varepsilon_1, (|x-y|/2)^{1/\rho}]} \quad \text{and} \quad 
		\{ H_\alpha(t\,,y) \}_{t \in [\varepsilon_1, (|x-y|/2)^{1/\rho}]}.
	\]
	By the Gaussian correlation inequality of Royen \cite{Royen},
	\begin{align*}
		\P\left\{ \sup_{t \in [\varepsilon_1, (|x-y|/2)^{1/\rho}]} 
		\frac{|H_\alpha(t\,,z)|}{t^{1/4}} \le \theta \right\} \le  
		\frac{\P\left\{ \sup_{t \in [\varepsilon_1, \varepsilon_2]} 
		\frac{|H_\alpha(t\,,z)|}{t^{1/4}} \le \theta \right\}}{%
		\P\left\{ \sup_{t \in [(|x-y|/2)^{1/\rho}, \varepsilon_2]} 
		\frac{|H_\alpha(t\,,z)|}{t^{1/4}} \le \theta \right\}},
	\end{align*}
	for $z \in \{x\,,y\}$.
	Choose and fix $\delta>0$.
	Then, the above, together with \eqref{varphi} and \eqref{Frostman:meas}, implies that
	\begin{align*}
		I_2 & \le \iint_{\substack{x,y \in K: \\ 2\varepsilon_1^\rho 
			< |x-y| \le 2\varepsilon_2^\rho}} \frac{\mu(\d x)\, \mu(\d y)}{%
			\P(\mathcal{E}((|x-y|/2)^{1/\rho}, \varepsilon_2 \,, x) ) 
			\P( \mathcal{E}((|x-y|/2)^{1/\rho}, \varepsilon_2\,,y) )}\\
		& \le \iint_{\substack{x,y \in K: \\ 2\varepsilon_1^\rho < |x-y| 
			\le 2\varepsilon_2^\rho}}  \left( \frac{\varepsilon_2}{%
			(|x-y|/2)^{1/\rho}} \right)^{2\lambda(\theta)+2\delta} 
			\mu(\d x) \, \mu(\d y)
			\lesssim \varepsilon_2^{2\lambda(\theta)+2\delta} 
			\varepsilon_1^{-2\lambda(\theta)-2\delta} \varepsilon_2^{\rho \beta},
	\end{align*}
	as $\max\{\varepsilon_2\,,\varepsilon_1/\varepsilon_2\} \to 0^+$, 
	subject to $\varepsilon_2^q \le \varepsilon_1\le \varepsilon_2$.
	We now choose $q>1$ sufficiently close to 1 so that 
	\[
		\rho\beta - (q-1) (2\lambda(\theta)+2\delta)>0,
	\]
	$h_{m+1} = h_m^q$, and $\max\{h_m\,,h_{m+1}/h_m\} \to 0^+$.
	Thanks to the above estimates, we have
	\[
		\E(X_{h_{m+1},h_m}^2) 
		\le 1+ \e^{-q^m(q\rho\beta+(q+1)\lambda(\theta) + 
		\mathscr{o}(1))} + \e^{-q^m(\rho\beta - (q-1)(2\lambda(\theta)+2\delta))} 
		 \text{ as $m\to\infty$.}
	\]
	Therefore, we may use the Paley-Zygmund inequality, 
	the fact that $\E(X_{h_{m+1},h_m})=1$, and the preceding to deduce that
	\begin{align*}
		&\P\{X_{h_{m+1},h_m}=0\}
		= 1- \P\{X_{h_{m+1},h_m}>0\}
			\le 1-\frac{[\E(X_{h_{m+1},h_m})]^2}{\E(X_{h_{m+1},h_m}^2)}\\
		&\qquad \le \e^{-q^m(q\rho\beta+(q+1)\lambda(\theta) + \mathscr{o}(1))} 
			+ \e^{-q^m(\rho\beta - (q-1)(2\lambda(\theta)+2\delta))} 
		\qquad \text{as $m\to\infty$.}
	\end{align*}
	It follows that $\sum_{m=1}^\infty \P\{X_{h_{m+1},h_m}=0\}<\infty$
	and the Borel-Cantelli lemma yields
	\begin{align*}
		\P\left( \bigcup_{n \in \N} \bigcap_{m \ge n} 
		\left\{ \exists\, x \in K\,, \sup_{t \in [h_{m+1},h_m]}
		\frac{|H_\alpha(t\,,x)|}{t^{1/4}} \le \theta \right\} \right)
		=1.
	\end{align*}
	This proves that there is a $\P$-null set off which there exists 
	$\mathfrak{N}=\mathfrak{N}(\omega) \in \N$ such that $U_{N,\mathfrak{N}}(\theta) \ne \varnothing$
	for all $N>\mathfrak{N}$, where $U_{N,n}(\theta)$ denotes the random subset of $K$
	that is defined as
	\begin{align*}
		U_{N,n}(\theta)  =  \left\{ x \in K : \sup_{t \in 
		[h_N,h_n]}\frac{|H_{\alpha}(t\,,x)|}{t^{1/4}}\le \theta \right\}
		\qquad\forall n\in\N,\ N>n.
	\end{align*}
	By the continuity of $H_{\alpha}$, every $U_{N,n}(\theta)$ 
	is a closed subset of the compact set $K$, and $N\mapsto U_{N,n}(\theta)$
	is decreasing with respect to set inclusion. Consequently,
	$\cap_{N \in \N: N>n} U_{N,n}(\theta)\neq\varnothing$ (Cantor's  theorem).
	This, together with Proposition \ref{pr:H-H_alpha:sup}, implies
	that there exists $x^*\in \cup_{n\in\N}
	\cap_{N \in \N: N>n} U_{N,n}(\theta)$ such that
	\[
		\adjustlimits\inf_{x \in K} \limsup_{t\to0^+} \frac{|H(t\,,x)|}{t^{1/4}}
		\le \limsup_{t\to0^+} \frac{%
		|H_{\alpha}(t\,,x^*)|}{t^{1/4}} \le \theta\quad\text{a.s.}
	\]
	Recall that $\lambda: (0\,,\infty) \to (0\,,\infty)$ is strictly decreasing and $\lambda(\theta) \to 0$ as $\theta \to \infty$ 
	(see Proposition \ref{pr:BC:H}).
	Therefore, we may let $\beta \uparrow \dimh K$, $\alpha \uparrow 1$ 
	and $\delta \downarrow 0$, all along rationals, in order to see that for every 
	$\theta>0$ that satisfies $\lambda(\theta) < \tfrac12 \dimh K$,
	\[
		\adjustlimits\inf_{x \in K} \limsup_{t\to0^+} \frac{|H(t\,,x)|}{t^{1/4}} \le \theta
		\quad\text{a.s.}
	\]
	Finally, we may take infimum over all such $\theta>0$ to 
	conclude that there is a $\P$-null set off which
	\[
		\adjustlimits\inf_{x \in K} \limsup_{t\to0^+} 
		\frac{|H(t\,,x)|}{t^{1/4}} \le \inf\{\theta>0 : \lambda(\theta)<\tfrac12 \dimh K\}
		=\lambda^{-1}\left(\tfrac12\dimh K\right)\quad\text{a.s.}
	\]
	This completes the proof of Theorem \ref{th:H}.
\end{proof} 

\section{Proof of Theorem \ref{th:BC:u}, Part (3)}
We will prove the slightly weaker result that
\begin{equation}\label{goal:dim}
	\P\left\{ \dimh \fS(\theta) = 
	1-2\lambda(\theta)\right\} =1 \qquad\forall \theta\ge\theta_c;
\end{equation}
see \eqref{theta_c}.
In fact, it follows from the above that
\[
	\P\left\{ \dimh \fS(\theta) = 1-2\lambda(\theta)
	\quad \forall \theta\in[\theta_c\,,\infty)\cap\Q_+\right\} =1,
\]
whence it follows from monotonicity that, off a single null set,
\[
	1-2\sup_{r\in\Q\cap[0,\theta]}
	\lambda(r) \le \dimh \fS(\theta) \le 1-2\inf_{r\in\Q\cap[\theta,\infty]}
	\lambda(r) 
	\qquad\forall \theta\ge\theta_c.
\]
The quantities on the left and the right of the above are both
equal to $\lambda(\theta)$, thanks to the continuity
of $\lambda$ (Proposition \ref{pr:BC:H}). Therefore, it remains to prove \eqref{goal:dim}.

The proof of \eqref{goal:dim}
hinges on a codimension argument, as was first introduced by
Taylor \cite[Theorem 4]{Taylor1966}. 
For every $\alpha\in(0\,,1)$ let $\{X_\alpha(t)\}_{t\ge0}$ denote
a symmetric $\alpha$-stable L\'evy process that is independent of the
random field $u$ in \eqref{u}.

Define $R(\alpha)$ to be the closure of the random set
$X_\alpha[0\,,1]$. Every $R(\alpha)$ is a random compact set since
$X_\alpha$ has c\`adl\`ag paths.
One can deduce from the potential theory of Hunt \cite{Hunt12,Hunt3}
that,  for every nonrandom Borel
(even analytic) set $L\subset\R$,
$\P\{ R(\alpha)\cap L\neq\varnothing\}>0$ if and only if
there exists a probability measure $\mu$ on $L$ such that
$\iint |a-b|^{-1+\alpha}\,\mu(\d a)\,\mu(\d b)<\infty$.
This and Frostman's theorem -- see Falconer \cite{Falconer} -- together yield
\begin{equation}\label{hit}
	\P\left\{ R(\alpha)\cap L\neq \varnothing\right\}
	\begin{cases}
		>0&\text{if $\dimh L>1-\alpha$,}\\
		=0&\text{if $\dimh L<1-\alpha$}.
	\end{cases}
\end{equation}
This is another way to say that  the ``codimension'' of the random set
$R(\alpha)$ is $1-\alpha$.

It is also a well-known fact that
\begin{equation}\label{s:hit}
	\dimh R(\alpha) = \dimm R(\alpha) = \alpha\quad\text{a.s.}
\end{equation}
The announced formula $\dimh R(\alpha) = \alpha$ follows from
Theorem 4.2 of Blumenthal and Getoor
\cite{BlumenthalGetoor}, and the formula $\dimm R(\alpha) = \alpha$
from Corollary 1.2 and Eq.~(1.5) of \cite{KSX12}.
Now we apply a codimension argument as follows: We apply
Part (2) of Theorem \ref{th:BC:u},
by first conditioning on $X_\alpha$, and then appeal to \eqref{s:hit} in order to see that
\[
	\P\left\{ \fS(\theta) \cap R(\alpha)\neq\varnothing \right\} =
	\begin{cases}
		1&\text{if $\lambda(\theta) < \alpha/2$},\\
		0&\text{if $\lambda(\theta) > \alpha/2$},
	\end{cases}
\]
for every $\theta>0$ and $\alpha\in(0\,,1)$. At the same time, 
\eqref{hit} also implies that, for every $\alpha\in(0\,,1)$ fixed,
\begin{align*}
	&\P\left( \fS(\theta) \cap R(\alpha)\neq\varnothing 
		\mid u\, \right)>0~\text{a.s.~on 
		$\left\{\omega\in\Omega:\ \dimh \fS(\theta)(\omega)>1-\alpha\right\}$},\\
	&\P\left(\fS(\theta) \cap R(\alpha)\neq\varnothing 
		\mid u\, \right)=0~\text{a.s.~on
		$\left\{\omega\in\Omega:\ \dimh \fS(\theta)(\omega)<1-\alpha\right\}$}.
\end{align*}
Together, the preceding two displays yield \eqref{goal:dim},
whence also Theorem \ref{th:BC:u}(3).
\qed

\bibliographystyle{plain}
\bibliography{SlowGrowth}

\begin{thebibliography}{10}

\bibitem{ACQ}
Gideon Amir, Ivan Corwin, and Jeremy Quastel.
\newblock Probability distribution of the free energy of the continuum directed
  random polymer in {$1+1$} dimensions.
\newblock {\em Comm. Pure Appl. Math.}, 64(4):466--537, 2011.

\bibitem{Anderson}
T.~W. Anderson.
\newblock The integral of a symmetric unimodal function over a symmetric convex
  set and some probability inequalities.
\newblock {\em Proc. Amer. Math. Soc.}, 6:170--176, 1955.

\bibitem{AthreyaJosephMueller}
Siva Athreya, Mathew Joseph, and Carl Mueller.
\newblock Small ball probabilities and a support theorem for the stochastic
  heat equation.
\newblock {\em Ann. Probab.}, 49(5):2548--2572, 2021.

\bibitem{BassBurdzy}
Richard~F. Bass and Krzysztof Burdzy.
\newblock A critical case for {B}rownian slow points.
\newblock {\em Probab. Theory Related Fields}, 105(1):85--108, 1996.

\bibitem{BlumenthalGetoor}
R.~M. Blumenthal and R.~K. Getoor.
\newblock Some theorems on stable processes.
\newblock {\em Trans. Amer. Math. Soc.}, 95:263--273, 1960.

\bibitem{Chen}
Jiaming Chen.
\newblock Small ball probabilities for the fractional stochastic heat equation
  driven by a colored noise.
\newblock {\em Electron. J. Probab.}, 30:Paper No. 35, 31, 2025.

\bibitem{Dalang1999}
Robert~C. Dalang.
\newblock Extending the martingale measure stochastic integral with
  applications to spatially homogeneous s.p.d.e.'s.
\newblock {\em Electron. J. Probab.}, 4:no. 6, 29, 1999.

\bibitem{DKN07}
Robert~C. Dalang, Davar Khoshnevisan, and Eulalia Nualart.
\newblock Hitting probabilities for systems of non-linear stochastic heat
  equations with additive noise.
\newblock {\em ALEA Lat. Am. J. Probab. Math. Stat.}, 3:231--271, 2007.

\bibitem{DS}
Robert~C. Dalang and Marta Sanz-Sol\'e.
\newblock {\em Stochastic partial differential equations, space-time white
  noise and random fields}.
\newblock Springer Monographs in Mathematics. Springer, Cham, [2026] \copyright
  2026.

\bibitem{Davis1983}
Burgess Davis.
\newblock On {B}rownian slow points.
\newblock {\em Z. Wahrsch. Verw. Gebiete}, 64(3):359--367, 1983.

\bibitem{DavisPerkins}
Burgess Davis and Edwin Perkins.
\newblock Brownian slow points: the critical case.
\newblock {\em Ann. Probab.}, 13(3):779--803, 1985.

\bibitem{Dvoretzky}
Aryeh Dvoretzky.
\newblock On the oscillation of the {B}rownian motion process.
\newblock {\em Israel J. Math.}, 1:212--214, 1963.

\bibitem{Falconer}
Kenneth Falconer.
\newblock {\em Fractal {G}eometry. {M}athematical {F}oundations and
  {A}pplications}.
\newblock John Wiley \& Sons, Ltd., Chichester, third edition, 2014.

\bibitem{FoondunJosephKim}
Mohammud Foondun, Mathew Joseph, and Kunwoo Kim.
\newblock Small ball probability estimates for the {H}\"older semi-norm of the
  stochastic heat equation.
\newblock {\em Probab. Theory Related Fields}, 185(1-2):553--613, 2023.

\bibitem{GRR}
A.~M. Garsia, E.~Rodemich, and H.~Rumsey, Jr.
\newblock A real variable lemma and the continuity of paths of some {G}aussian
  processes.
\newblock {\em Indiana Univ. Math. J.}, 20:565--578, 1970/71.

\bibitem{GarsiaRR}
A.~M. Garsia, E.~Rodemich, and H.~Rumsey, Jr.
\newblock A real variable lemma and the continuity of paths of some {G}aussian
  processes.
\newblock {\em Indiana Univ. Math. J.}, 20:565--578, 1970/71.

\bibitem{GreenwoodPerkins}
Priscilla Greenwood and Edwin Perkins.
\newblock A conditioned limit theorem for random walk and {B}rownian local time
  on square root boundaries.
\newblock {\em Ann. Probab.}, 11(2):227--261, 1983.

\bibitem{GuoSongWangXiao}
Yuhui Guo, Jian Song, Ran Wang, and Yimin Xiao.
\newblock Sample path properties and small ball probabilities for stochastic
  fractional diffusion equations.
\newblock {\em J. Differential Equations}, 446:Paper No. 113604, 56, 2025.

\bibitem{HairerPardoux}
Martin Hairer and \'Etienne Pardoux.
\newblock A {W}ong-{Z}akai theorem for stochastic {PDE}s.
\newblock {\em J. Math. Soc. Japan}, 67(4):1551--1604, 2015.

\bibitem{HuLee}
Jingwu Hu and Cheuk~Yin Lee.
\newblock On the spatio-temporal increments of nonlinear parabolic {SPDEs} and
  the open {KPZ} equation, 2025.
\newblock To appear in Ann. Appl. Probab. E-print available at
  \url{https://arxiv.org/abs/2508.05032}.

\bibitem{HuangKh}
Jingyu Huang and Davar Khoshnevisan.
\newblock On the multifractal local behavior of parabolic stochastic {PDE}s.
\newblock {\em Electron. Commun. Probab.}, 22:Paper No. 49, 11, 2017.

\bibitem{Hunt12}
G.~A. Hunt.
\newblock Markoff processes and potentials. {I}, {II}.
\newblock {\em Illinois J. Math.}, 1:44--93, 316--369, 1957.

\bibitem{Hunt3}
G.~A. Hunt.
\newblock Markoff processes and potentials. {III}.
\newblock {\em Illinois J. Math.}, 2:151--213, 1958.

\bibitem{Kahane1974}
Jean-Pierre Kahane.
\newblock Sur l'irr\'egularit\'e{} locale du mouvement brownien.
\newblock {\em C. R. Acad. Sci. Paris S\'er. A}, 278:331--333, 1974.

\bibitem{Kahane1976}
Jean-Pierre Kahane.
\newblock Sur les z\'eros et les instants de ralentissement du mouvement
  brownien.
\newblock {\em C. R. Acad. Sci. Paris S\'er. A-B}, 282(8):Aii, A431--A433,
  1976.

\bibitem{Kahane1983}
Jean-Pierre Kahane.
\newblock Slow points of {G}aussian processes.
\newblock In {\em Conference on {H}armonic {A}nalysis in {H}onor of {A}ntoni
  {Z}ygmund, {V}ol. {I}, {II} ({C}hicago, {I}ll., 1981)}, Wadsworth Math. Ser.,
  pages 67--83. Wadsworth, Belmont, CA, 1983.

\bibitem{KSX12}
D.~Khoshnevisan, R.~L. Schilling, and Y.~Xiao.
\newblock Packing dimension profiles and {L}\'evy processes.
\newblock {\em Bull. Lond. Math. Soc.}, 44(5):931--943, 2012.

\bibitem{CBMS}
Davar Khoshnevisan.
\newblock {\em Analysis of {S}tochastic {P}artial {D}ifferential {E}quations},
  volume 119 of {\em CBMS Regional Conference Series in Mathematics}.
\newblock Conference Board of the Mathematical Sciences, Washington, DC; by the
  American Mathematical Society, Providence, RI, 2014.

\bibitem{KKimMueller}
Davar Khoshnevisan, Kunwoo Kim, and Carl Mueller.
\newblock Small-ball constants, and exceptional flat points of {SPDE}s.
\newblock {\em Electron. J. Probab.}, 29:Paper No. 180, 31, 2024.

\bibitem{KL}
Davar Khoshnevisan and Cheuk~Yin Lee.
\newblock On the passage times of self-similar {G}aussian processes on curved
  boundaries.
\newblock {\em Bernoulli}, 32(3):2098--2122, 2026.

\bibitem{KSXZ}
Davar Khoshnevisan, Jason Swanson, Yimin Xiao, and Liang Zhang.
\newblock Weak existence of a solution to a differential equation driven by a
  very rough f{B}m, 2013.
\newblock unpublished manuscript. Available electronically on
  \url{https://arxiv.org/abs/1309.3613}.

\bibitem{Ledoux}
Michel Ledoux.
\newblock {\em The {C}oncentration of {M}easure {P}henomenon}, volume~89 of
  {\em Mathematical Surveys and Monographs}.
\newblock American Mathematical Society, Providence, RI, 2001.

\bibitem{Martin}
A.~Martin.
\newblock Small ball asymptotics for the stochastic wave equation.
\newblock {\em J. Theoret. Probab.}, 17(3):693--703, 2004.

\bibitem{Perkins1983}
Edwin Perkins.
\newblock On the {H}ausdorff dimension of the {B}rownian slow points.
\newblock {\em Z. Wahrsch. Verw. Gebiete}, 64(3):369--399, 1983.

\bibitem{Royen}
Thomas Royen.
\newblock A simple proof of the {G}aussian correlation conjecture extended to
  some multivariate gamma distributions.
\newblock {\em Far East J. Theor. Stat.}, 48(2):139--145, 2014.

\bibitem{Taylor1966}
S.~J. Taylor.
\newblock Multiple points for the sample paths of the symmetric stable process.
\newblock {\em Z. Wahrscheinlichkeitstheorie und Verw. Gebiete}, 5:247--264,
  1966.

\bibitem{Walsh}
John~B. Walsh.
\newblock An {I}ntroduction to {S}tochastic {P}artial {D}ifferential
  {E}quations.
\newblock In {\em \'Ecole d'\'et\'e{} de {P}robabilit\'es de {S}aint-{F}lour,
  {XIV}---1984}, volume 1180 of {\em Lecture Notes in Math.}, pages 265--439.
  Springer, Berlin, 1986.

\end{thebibliography}

\end{document}